\documentclass[a4paper,11pt,reqno]{article}
\usepackage{amsmath, amsfonts, amssymb, amsthm}
\input{epsf.sty}
\usepackage{mathrsfs}
\usepackage{amsmath}
\usepackage{amssymb,color}
\usepackage{graphicx}
\usepackage{subfigure}
\usepackage{amssymb}

\usepackage[colorlinks=true]{hyperref}
\hypersetup{urlcolor=blue, citecolor=blue}

\newtheorem{proposition}{Proposition}

\def \R {{\mathbb{R}}}
\numberwithin{equation}{section}

\begin{document}

\title{Stability of peakons for the generalized modified Camassa-Holm equation}

\author{
Zihua $\mbox{Guo}^{1}$ \footnote{E-mail: zihua.guo@monash.edu}
\quad  Xiaochuan $\mbox{Liu}^{2}$ \footnote{E-mail: liuxiaochuan@mail.xjtu.edu.cn}
\quad  Xingxing $\mbox{Liu}^{1,3}$ \footnote{Corresponding author. E-mail: liuxxmaths@cumt.edu.cn}
and 
\quad Changzheng $\mbox{Qu}^{4}$ \footnote{E-mail: quchangzheng@nbu.edu.cn}
\\
$^1\mbox{School}$ of Mathematical Sciences, Monash University,
\\
Melbourne, VIC 3800, Australia
\\
$^2\mbox{School}$ of Mathematics and Statistics, Xi'an Jiaotong University,
\\
Xi'an, Shaanxi 710049, China
\\
$^3\mbox{School}$ of mathematics, China University of Mining and Technology,
\\Xuzhou, Jiangsu 221116, China
\\
$^4\mbox{Center}$ for Nonlinear Studies and Department of Mathematics, Ningbo University,
\\
Ningbo, Zhejiang 315211, China
}

\date{}
\maketitle

\begin{abstract}
In this paper, we study orbital stability of peakons for the generalized modified Camassa-Holm (gmCH) equation, which is a natural higher-order generalization of
the modified Camassa-Holm (mCH) equation, and admits Hamiltonian form and single peakons. We first show that the single peakon is the usual weak solution of the PDEs. Some sign invariant properties and conserved densities are presented. Next, by constructing the corresponding auxiliary function $h(t,\,x)$ and establishing a delicate polynomial inequality relating to the two conserved densities with the maximal value of approximate solutions, the orbital stability of single peakon of the gmCH equation is verified. We introduce a new approach to prove the key inequality, which is different from that used for the mCH equation. This extends the result on the stability of peakons for the mCH equation (Comm. Math. Phys., 322:967-997, 2013) successfully to the higher-order case, and is helpful to understand how higher-order nonlinearities affect the dispersion dynamics.
\\


\end{abstract}

\section{Introduction}
\par

In the past two decades, the Camassa-Holm (CH) equation
\begin{equation*}
y_t+uy_x+2u_x y=0, \ \ \ y=u-u_{xx},
\end{equation*}
attracted much attention among the communities of the nonlinear integrable equations and the PDEs.
In 1993, Camassa and Holm \cite{C-H} obtained the CH equation
by approximating directly in the Hamiltonian for Euler's equations in the shallow water regime.
It can model the unidirectional propagation of shallow water waves over a flat bottom \cite{C-H,Constantin-L,Johnson},
with $u(t,x)$ standing for the fluid velocity at time $t\geq 0$ in the spatial $x\in \R$ direction.
Actually, the CH equation was initially introduced in 1981 by Fuchssteiner and
Fokas \cite{F-Fokas} as a bi-Hamiltonian generalization of the KdV equation.
Dai \cite{Dai} derived it as an equation describing the propagation of axially symmetric waves in hyperelastic rods.
The local well-posedenss of the CH equation in $H^s$ $(s>3/2)$ was shown in \cite{Li-Olver}, while ill-posedness were proved in \cite{Byers}
and \cite{G-L-M-Y} for $s<3/2$ and $s=3/2$, respectively.
The CH equation, like the classical KdV equation, admits bi-Hamiltonian structure \cite{F-Fokas}  and is completely
integrable \cite{Beals-S-S,Constantin1}. However, the CH equation possesses several nontrivial properties. Concerning about the global solutions, it is well-known that any smooth solutions of the KdV equation are global, and the global strong solutions of the CH equation may describe wave breaking \cite{Constantin0,C-E1,C-E2,C-E3}, i.e., the wave profile remains bounded, but its slope becomes unbounded in finite time \cite{Whitham}.

Another remarkable property of the CH equation is the presence of peaked soliton solutions \cite{C-H-H}, called \emph{peakons}.
They are given by
\begin{eqnarray*}
u(t,x)=\varphi_c(x-c\,t)=\varphi(x-c\,t)=ce^{-|x-c\,t|},\ \ \ c\in \R,
\end{eqnarray*}
which are solitons, retaining their shapes and speeds after interacting with other peakons \cite{Beals-S-S1}.
It is worth pointing out that the feature of peakons that their profile is smooth except for a peak at its crest,
is analogous to that of the waves of greatest height, i.e., traveling waves of largest possible amplitude which are
solutions to be governing equations for water waves \cite{Constantin2,C-E4,Toland}. It is well-known that the stability
of solitary waves is one of the fundamental questions for nonlinear wave equations.
However, the stability of peakons to the CH equation seems not to enter the general framework developed for
instance in \cite{Benjamin,Grillakis-S-S}. In an intriguing paper \cite{Constantin-S}, Constantin and Strauss provided a direct proof to orbital stability of peakons of the CH equation by employing its conserved densities and the specific structure of the peakons. In 2009, El Dika and Molinet \cite{El Dika-M}
proved orbital stability of multi-peakons by combining the proof of stability of single peakons with
a property of almost monotonicity of the localized energy. Furthermore, they also investigated the orbital stability of ordered trains of anti-peakons and peakons  \cite{El Dika-M1}.

Recently, the great interest in the CH equation has inspired the
search for various CH-type equations with cubic or higher-order nonlinearities. One of the most
concerned is the following modified CH (mCH) equation
\begin{equation*}
y_t+\big((u^2-u_x^2)y\big)_x=0,\ \ \ y=u-u_{xx},
\end{equation*}
which was derived by Fuchssteiner \cite{Fuchssteiner} and Olver and Rosenau \cite{O-R}
by employing the tri-Hamiltonian duality approach to the bi-Hamiltonian representation of the modified KdV equation. Subsequently, the integrability and structure of solutions to the mCH equation was discussed by Qiao \cite{Qiao1}. The mCH equation is also called FORQ equation in some literature \cite{H-Man,Yang-L-Z}.

The mCH equation is completely integrable \cite{O-R}. It has a bi-Hamiltonian structure and also admits a Lax pair \cite{Qiao3},
and hence may be solved by the inverse scattering transform method. Compared with the CH equation, the mCH equation
admits not only peakons, but also possesses cusp solitons (cuspons) and weak kink solutions
 ($u,u_x,u_t$ are continuous, but $u_{xx}$
has a jump at its peak point) \cite{Qiao3,Qiao4}. It has also
significant differences from the CH equation on the dynamics of
the multi-peakons and peakon-kink solutions \cite{G-L-O-Q, Qiao2, Qiao4}.
Fu et al., \cite{Fu} studied the Cauchy problem of the mCH equation in Besov spaces
and the blow-up scenario. The non-uniform dependence on the initial data was established in \cite{H-Man}.
Gui et al., \cite{G-L-O-Q} considered the formulation of singularities of solutions and showed
that some solutions with certain initial data would blow up in finite
time. Then the blow-up phenomena were systematically investigated in \cite{C-L-Q-Z,L-O-Q-Z}.
The mCH equation admits single peakon of the form \cite{G-L-O-Q}
\begin{eqnarray*}
u(t,x)=\varphi_c(x-c\,t)=\sqrt{\frac{3\,c}{2}}e^{-|x-c\,t|},\quad c>0.
\end{eqnarray*}
Later, inspired by \cite{Constantin-S} and \cite{El Dika-M}, the orbital stability of single peakon and the train
of peakons for the mCH equation were proved in \cite{Qu} and \cite{L-L-Q}, respectively.

In this paper, we consider the following generalized modified CH (gmCH) equation proposed in \cite{R-A}:
\begin{equation}\label{1.1}
y_t+\big((u^2-u_x^2)^ny\big)_x=0,\ \ \ y=u-u_{xx},
\end{equation}
where $n\geq 2$ is a positive integer. When $n=1$, equation (\ref{1.1}) becomes the mCH equation.
Motivated by the fact that the generalized KdV equation, $v_t-v^pv_x-v_{xxx}=0 \; (p\geq 1)$,
shared one of the KdV's two Hamiltonian structure and reduced to the KdV equation when $p=1,$
Recio and Anco \cite{R-A} derived the gmCH equation (\ref{1.1}) in the classification of the family of nonlinear dispersive wave equations involving two arbitrary functions when they admit peaked solitary waves.
Very recently, the local well-posedness in Besov spaces for the Cauchy problem associated to equation (\ref{1.1}), and blow-up mechanism have been discussed
in \cite{Yang-L-Z}. On the other hand, one of the main interesting features of the CH equation (with quadratic nonlinearity) and the mCH equation (with cubic nonlinearity) is the existence of stable peakons. Notice that, the gmCH equation (\ref{1.1}) with higher-order nonlinearity can also have single peakon and multi-peakons \cite{R-A}. Thus it is of great interest to study whether the single peakon of equation (\ref{1.1}) is orbitally stable as the CH and mCH equations.

To pursue the above goal, our approach to prove the orbital stability of peakons for equation (\ref{1.1}) is adapted from
the ingenious idea of Constantin and Strauss in \cite{Constantin-S}, and evolved for the mCH equation \cite{Qu}. However, due to the generalized
form with higher-order nonlinearity, the analysis of stability result for the gmCH equation (\ref{1.1}) is
more subtle than that for the CH and mCH equations. On the one hand, owing to the same conservation law
$E(u)$ as the CH and mCH equations, we naturally expect the orbital stability of peakons for equation (\ref{1.1})
in the sense of the energy space $H^1(\R)$-norm. Note that equation \eqref{1.1} can be written as the Hamitonian form \cite{R-A}, where the Hamiltonian conserved density 
\begin{eqnarray*}
H_{gmCH}(u)=\frac 1{2(n+1)}\int_{\R}u(u^2-u^2_x)^n ydx,\ \ \ y=u-u_{xx},
\end{eqnarray*}
is more complicated than that of the mCH equation. For our purpose, we derive an equivalent form of $F(u)$ only involving the terms of $u$ and $u_x$. Note that the two useful conservation laws $E(u)$ and $F(u)$ play a crucial role in our approach.

On the other hand, in view of \cite{Constantin-S,Qu}, we see that the first step in the proof of orbital stability
is to construct two auxiliary functionals $g$ and $h$, so that we can establish an inequality relating to the Hamiltonian conserved densities $E(u)$ and $F(u)$ as well as the maximal value of approximate solutions $u$ (see Lemma 3.2). There are two difficulties encountered by equation (\ref{1.1}) in establishing the polynomial inequality.
First, although we can define the same functional $g$ as in the CH and mCH equations because of the same conserved quantity $E(u)$, here we need more subtle analysis to construct the other functional $h$. Fortunately, based on the observation that the functionals $g$ and $h$ are required to vanish at the peakons, we can obtain a polynomial of degree $2n$
as $h$ from the equivalent form of $F(u)$. The second difficulty, arising from the derivation of the required polynomial inequality, is to show
$h\leq \frac{2-c_1}{2} M^{2n}, M\triangleq \max_{x\in \R}u(x)$ (see (\ref{3.10})),
which is substantially different from the cases for the CH equation ($h=u\leq M$) and the mCH equation ($h=u^2\mp\frac{2}{3}uu_x-\frac{1}{3}u^2_x\leq \frac{4}{3} M^2$). This new difficulty is caused by the higher-order nonlinear structure and higher-order conservation laws. To tackle it, by carefully analyzing the properties of the coefficients $c_k,d_k$ of the function $h$, we equivalently convert the estimate of $h$ into the non-positivity of some function.
Then we complete the proof of it by using delicate analytic calculations (see Lemma 3.2).
In addition, we need the sign-invariant property of $y$, similar to the case of $n=1$ of equation (\ref{1.1}) in \cite{Qu},  to
guarantee that the solution $y$ is positive and $(1\pm \partial_x)u\geq 0$ when the initial value $y_0(x)\geq 0$. This enables us to prove the main inequality as follows
\begin{equation*}
\frac{n(2-c_1)}{n+1}M^{2n+2}-\frac{2-c_1}{2}M^{2n}E(u)+F(u)\leq 0,
\end{equation*}
where the constant $c_1$ is related to the auxiliary function $h$. Indeed, after overcoming the above difficulties, the strategy of proof of stability is rather routine. By expanding $E(u)$ around the peakon $\varphi_c$, the $H^1$-norm between the solution $u$ and the
peakon $\varphi_c(\cdot-\xi(t))$ is precisely controlled by the error term $|M-\max_{x\in \R}\varphi_c|$,
while the error term can just be estimated by analyzing the root structure of the obtained polynomial inequality.

The remainder of this paper is organized as follows. In Section 2, we show the existence of peaked solitons
which can be understood as weak solutions of equation (\ref{1.1}), and then present two conserved quantities which are crucial
in the proof of stability. In Section 3, the orbital stability of peaked solitons for equation (\ref{1.1}) is verified based on several useful lemmas. In Section 4, we end the paper with an appendix devoted to the proofs of the last third equality (\ref{3.8}), last inequality of (\ref{3.29}) and the conservation law $F(u)$.

$Notation.$ In the followings, we use the notation $A\lesssim B$ to denote the corresponding inequality $A\leq C B$  for some constant $C>0$.

\section{Preliminaries}

\newtheorem {remark2}{Remark}[section]
\newtheorem {proposition2}{Proposition}[section]
\newtheorem {definition2}{Definition}[section]
\newtheorem{theorem2}{Theorem}[section]
\newtheorem{lemma2}{Lemma}[section]

In this section, we first recall the local well-posedness result to the Cauchy problem of the gmCH equation (\ref{1.1}) (up to a slight modification), and properties for strong solutions. Then we show the existence of single peakon and present two useful conserved densities, which will be frequently used in the rest of the paper.
\begin{lemma2} \cite{Yang-L-Z}
Let $u_0(x)=u(0,x)\in H^s(\R)$ with $s>5/2$. Then there exists a time $T>0$ such that
the Cauchy problem of (\ref{1.1}) has a unique strong solution $u(t,x)\in C([0,T);H^s(\R))\cap
C^1([0,T);H^{s-1}(\R))$ and the map $u_0\mapsto u$ is continuous from a neighborhood of $u_0$ in $H^s(\R)$ into
$C([0,T);H^s(\R))\cap C^1([0,T);H^{s-1}(\R))$.
\end{lemma2}

The conservative property of the momentum density $y$ has been proved to play a crucial role in the study of well-posedness of global solutions and formation of singularities to the CH and mCH equations. For the gmCH equation, we consider the flow governed by the effective wave speed $(u^2-u_x^2)^n$
\begin{equation}\label{hmCH-flow}
\begin{cases}
\frac{dq(t,x)}{dt}= \big(u^2(t,q(t,x)) - u^2_x(t,q(t,x)) \big)^n, \;\; x \in \mathbb{R},\quad t\in [0,T).\\
q(0,x)=x,\;\; x \in \mathbb{R}.
\end{cases}
\end{equation}
Differentiating the equation in \eqref{hmCH-flow} with respect to $x$ leads to the following result.

\begin{proposition}\label{lem:conservation}
Suppose $u_0 \in H^s(\mathbb{R})$ with  $s>5/2$, and let $T>0$ be the maximal existence
time of the strong solution $u(t,x)\in C\left([0, T), H^{s}({\Bbb R})\right) \cap C^1 \left([0, T), H^{s-1}({\Bbb R})\right)$ to the Cauchy problem of \eqref{1.1}. Then \eqref{hmCH-flow} has a unique solution
$q \in C^1([0,T)\times \mathbb{R},\mathbb{R})$ such that $q(t,\cdot)$ is an increasing diffeomorphism over $\R$ with
\begin{equation}\label{deriv of q}
q_x(t,x)=\exp\left(2n\int_0^t ((u^2-u_x^2)^{n-1} u_xy )(\tau,q(\tau,x)) \,d\tau\right)>0
\end{equation}
for all $(t,x)\in [0,T)\times\mathbb{R}$. Furthermore, the momentum density $y=u - u_{xx}$ satisfies
\begin{equation}\label{cons_m}
y(t, q(t,x)) q_x(t,x) = y_0(x), \qquad (t,x)\in [0,T) \times \mathbb{R},
\end{equation}
which implies that the sign and zeros of $y$ are preserved under the flow.
\end{proposition}

The proof of this lemma is similar to that of Lemma 5.1 for the mCH equation in \cite{G-L-O-Q}, we omit it here. Using the above proposition, we can prove the following lemma, which is similar to Lemmas 2.8-2.9 in \cite{L-L}.
\begin{lemma2}
Assume $u_0(x)\in H^s(\R)$, $s>5/2.$ If $y_0(x)=(1-\partial^2_x)u_0(x)$ does not change sign, then $y(t,x)$
will not change sign for all $t\in [0,T)$. It follows that if $y_0\geq0$, then the corresponding solution $u(t,x)$
of equation (\ref{1.1}) is positive for $(t,x)\in[0,T)\times \R. $
Furthermore, if $y_0\geq0$, then the corresponding solution $u(t,x)$ of equation (\ref{1.1}) satisfies
\begin{eqnarray*}
(1\pm \partial_x)u(t,x)\geq 0, \ \quad \ \mbox{for}\ \quad \forall(t,x)\in[0,T)\times \R.
\end{eqnarray*}
\end{lemma2}
Now, with $y=u-u_{xx},$ we can rewrite the gmCH equation (\ref{1.1}) in the following more convenient form
\begin{eqnarray}\label{2.0}
&&u_t-u_{txx}+\sum_{k=0}^{n-1}(-1)^k\big(C^k_n+2nC_{n-1}^k\big)u^{2n-2k}u_x^{2k+1}+(-1)^nu_x^{2n+1}\nonumber\\
&&-\sum_{k=0}^{n-1}(-1)^k4nC_{n-1}^ku^{2n-2k-1}u_x^{2k+1}u_{xx}+
\sum_{k=0}^{n-1}(-1)^k2nC_{n-1}^ku^{2n-2k-2}u_x^{2k+1}u^2_{xx}
\nonumber\\
&&-\sum_{k=0}^{n}(-1)^kC_{n}^ku^{2n-2k}u_x^{2k}u_{xxx}=0.
\end{eqnarray}
Applying the operator $(1-\partial_x^2)^{-1}$ to both sides of the above equation, one gets
\begin{eqnarray}\label{2.1}
&&u_t+\left(u^{2n}-\sum_{k=1}^n\frac{(-1)^{k+1}}{2k+1}C^k_nu^{2n-2k}u_{x}^{2k}\right)u_x+
(1-\partial_x^2)^{-1}\partial_x\Bigg(\frac{2n}{2n+1}u^{2n+1}\nonumber\\
&+&\sum_{k=1}^n\frac{(-1)^{k-1}}{2k-1}C^k_nu^{2n-2k+1}u_{x}^{2k}\Bigg)
+(1-\partial_x^2)^{-1}\left(\sum_{k=1}^n\frac{(-1)^{k+1}}{2k+1}C^k_nu^{2n-2k}u_{x}^{2k+1}\right)=0.
\end{eqnarray}
Recall that $(1-\partial_x^2)^{-1}f=p\ast f$ for all $f\in L^2$, where $p(x)\triangleq e^{-|x|}/2, x\in\R $. We now apply (\ref{2.1}) to define weak solutions of equation (\ref{1.1}).

\begin{definition2}
Let $u_0\in W^{1,{2n+1}}(\R)$ be given. If $u(t,x)\in
L^{\infty}_{loc}([0,T);W_{loc}^{1,{2n+1}}(\R))$ and satisfies
\begin{eqnarray*}
&&\int_0^T\int_{\R}\Bigg[u\phi_t+\frac{1}{2n+1}u^{2n+1}\phi_x+\left(\sum_{k=1}^n\frac{(-1)^{k+1}}{2k+1}C^k_nu^{2n-2k}u_{x}^{2k+1}\right)
\phi\\
&+&
p\ast\left(\frac{2n}{2n+1}u^{2n+1}
+\sum_{k=1}^n\frac{(-1)^{k-1}}{2k-1}C^k_nu^{2n-2k+1}u_{x}^{2k}\right)\cdot\partial_x\phi\\
&-&p\ast
\left(\sum_{k=1}^n\frac{(-1)^{k+1}}{2k+1}C^k_nu^{2n-2k}u_{x}^{2k+1}\right)\cdot\phi\Bigg]dxdt +\int_{\R}u_0(x)\phi(0,x)dx=0,
\end{eqnarray*}
for any $\phi(t,x)\in C_c^{\infty}([0,T)\times\R),$  then
$u(t,x)$ is called a weak solution to equation (\ref{1.1}). If $u$ is a
weak solution on $[0,T)$ for every $T>0$, then it is called a global
weak solution.
\end{definition2}

\begin{remark2}
Since the Sobolev space $W_{loc}^{1,{2n+1}}(\R))$ can be embedded in the H\"{o}lder space $C^\alpha(\R)$
with $0\leq \alpha \leq \frac{2n}{2n+1}$, Definition 2.1 precludes the admissibility of discontinuous shock
waves as weak solutions.
\end{remark2}
Based on Definition 2.1, we show the existence of single peakon to equation \eqref{1.1}.
\begin{theorem2}
For any $c>0$, the peaked functions of the form
\begin{equation}\label{2.2}
\varphi_c(t,x)=a\,e^{-|x-ct|}, \quad \ \mbox{where}\quad c=\left(1-\sum_{k=1}^n\frac{(-1)^{k+1}}{2k+1}C^k_n\right)a^{2n},
\end{equation}is a global weak solution to the gmCH equation (\ref{1.1}) in the sense of Definition 2.1.
\end{theorem2}
\begin{proof}
For any test function $\phi(\cdot)\in C_c^{\infty}(\R)$, we infer
\begin{eqnarray*}
\int_{\R}e^{-|x|}\phi'(x)dx&=&\int_{-\infty}^0e^x\phi'(x)dx+\int_0^{\infty}e^{-x}\phi'(x)dx\\
&=&e^x\phi(x)\Big|^0_{-\infty}-\int_{-\infty}^0e^x\phi(x)dx+e^{-x}\phi(x)\Big|^{\infty}_{0}
+\int_0^{\infty}e^{-x}\phi(x)dx\\
&=&-\int_{-\infty}^0e^x\phi(x)dx+\int_0^{\infty}e^{-x}\phi(x)dx=\int_{\R} \mbox{sign} (x)e^{-|x|}\phi(x)dx,
\end{eqnarray*}
after integration by parts. Thus, for all $t\geq0$, we have
\begin{eqnarray}\label{2.3}
\partial_x\varphi_c(t,x)=-\mbox{sign}(x-ct)\varphi_c(t,x),
\end{eqnarray}
in the sense of distribution $\mathcal{S}'(\R)$.
Letting $\varphi_{0, c}(x)\triangleq \varphi_c(0,x)$, then we get
\begin{eqnarray}\label{2.4}
\lim \limits_{t\rightarrow
0^+}\|\varphi_c(t,\cdot)-\varphi_{0,c}(x)\|_{W^{1,\infty}}=0.
\end{eqnarray}
The same computation as in (\ref{2.3}), for all $t\geq0$, yields,
\begin{eqnarray}\label{2.5}
\partial_t\varphi_c(t,x)=c\ \mbox{sign}(x-ct)\varphi_c(t,x)\in L^\infty.
\end{eqnarray}
Combining (\ref{2.3})-(\ref{2.5}) with integration by parts, for
any test function $\phi(t,x)\in C_c^{\infty}([0,\infty)\times
\R)$, we obtain
\begin{eqnarray}\label{2.6}
&&\int_0^{+\infty}\int_{\R}\Big[\varphi_c\partial_t\phi+\frac{1}{2n+1}\varphi_c^{2n+1}\partial_x\phi
+\Big(\sum_{k=1}^n\frac{(-1)^{k+1}}{2k+1}C^k_n\varphi_c^{2n-2k}(\partial_x\varphi_c)^{2k+1}\Big)
\phi\Big]dxdt\nonumber\\
&&+\int_{\R}\varphi_{0, c}(x)\phi(0,x)dx\nonumber\\
&=&-\int_0^{+\infty}\int_{\R}\Big[\partial_t\varphi_c+\varphi_c^{2n}\partial_x\varphi_c-
\Big(\sum_{k=1}^n\frac{(-1)^{k+1}}{2k+1}C^k_n\varphi_c^{2n-2k}(\partial_x\varphi_c)^{2k+1}\Big)\Big]\phi dx dt \nonumber\\
&=&-\int_0^{+\infty}\int_{\R}\phi\
\mbox{sign}(x-ct)\varphi_c
\Big[c-\Big(1-\sum_{k=1}^n\frac{(-1)^{k+1}}{2k+1}C^k_n\Big)\varphi_c^{2n}\Big]dxdt.
\end{eqnarray}
Using the definition of $\varphi_c$ and the relation $c=\left(1-\sum\limits_{k=1}^n\frac{(-1)^{k+1}}{2k+1}C^k_n\right)a^{2n}$, we deduce that for $x>ct$,
\begin{eqnarray}\label{T1}
&&\mbox{sign}(x-ct)\varphi_c
\Big[c-\Big(1-\sum_{k=1}^n\frac{(-1)^{k+1}}{2k+1}C^k_n\Big)\varphi_c^{2n}\Big]\nonumber\\
&=&a^{2n+1}\Big(1-\sum_{k=1}^n\frac{(-1)^{k+1}}{2k+1}C^k_n\Big)
\big(e^{ct-x}-e^{(2n+1)(ct-x)}\big),
\end{eqnarray}
and for $x\leq ct$,
\begin{eqnarray}\label{T2}
&&\mbox{sign}(x-ct)\varphi_c
\Big[c-\Big(1-\sum_{k=1}^n\frac{(-1)^{k+1}}{2k+1}C^k_n\Big)\varphi_c^{2n}\Big]\nonumber\\
&=&-\Big(1-\sum_{k=1}^n\frac{(-1)^{k+1}}{2k+1}C^k_n\Big)
a^{2n+1}\big(e^{x-ct}-e^{(2n+1)(x-ct)}\big).
\end{eqnarray}
On the other hand, by Definition 2.1, we derive
\begin{eqnarray}\label{2.7}
\begin{aligned}
&\int_0^{+\infty}\int_{\R}\Big[(1-\partial_x^2)^{-1}\Big(\frac{2n}{2n+1}\varphi_c^{2n+1}+
\sum_{k=1}^n\frac{(-1)^{k-1}}{2k-1}C^k_n\varphi_c^{2n-2k+1}(\partial_x\varphi_c)^{2k}\Big)\partial_x\phi\\
&\qquad\qquad\quad-(1-\partial_x^2)^{-1}\Big(\sum_{k=1}^n\frac{(-1)^{k+1}}{2k+1}C^k_n\varphi_c^{2n-2k}(\partial_x\varphi_c)^{2k+1}\Big)\phi\Big]dx
dt\\
&=-\int_0^{+\infty}\int_{\R}\Big[\phi\cdot\partial_x p\ast\big(\sum_{k=1}^n\frac{(-1)^{k-1}}{2k-1}C^k_n\varphi_c^{2n-2k+1}(\partial_x\varphi_c)^{2k}\big)\\
&\qquad\qquad\qquad\quad+\phi\cdot p\ast \big(2n \varphi_c^{2n}\partial_x\varphi_c+\sum_{k=1}^n\frac{(-1)^{k+1}}{2k+1}C^k_n\varphi_c^{2n-2k}(\partial_x\varphi_c)^{2k+1}\big)\Big]dxdt.
\end{aligned}
\end{eqnarray}
We calculate from (\ref{2.3}) that,
\begin{eqnarray*}
&&2n \varphi_c^{2n}\partial_x\varphi_c+\sum_{k=1}^n\frac{(-1)^{k+1}}{2k+1}C^k_n\varphi_c^{2n-2k}(\partial_x\varphi_c)^{2k+1}\nonumber\\
&=&2n\varphi_c^{2n}\big(-\mbox{sign}(x-ct)\varphi_c\big)+\sum_{k=1}^n\frac{(-1)^{k+1}}{2k+1}C^k_n\varphi_c^{2n-2k}
\big(-\mbox{sign}(x-ct)\varphi_c\big)^{2k+1}\nonumber\\
&=&\Big[2n+\sum_{k=1}^n\frac{(-1)^{k+1}}{2k+1}C^k_n\Big](-\mbox{sign}(x-ct))\varphi_c^{2n+1}\nonumber=\frac{1}{2n+1}\Big[2n+\sum_{k=1}^n\frac{(-1)^{k+1}}{2k+1}C^k_n\Big]\partial_x(\varphi_c^{2n+1}),
\end{eqnarray*}
which together with (\ref{2.7}) leads to
\begin{eqnarray}\label{2.8}
\begin{aligned}
&\int_0^{+\infty}\int_{\R}\Big[(1-\partial_x^2)^{-1}\Big(\frac{2n}{2n+1}\varphi_c^{2n+1}+
\sum_{k=1}^n\frac{(-1)^{k-1}}{2k-1}C^k_n\varphi_c^{2n-2k+1}(\partial_x\varphi_c)^{2k}\Big)\partial_x\phi\\
&\qquad\qquad\quad-(1-\partial_x^2)^{-1}\Big(\sum_{k=1}^n\frac{(-1)^{k+1}}{2k+1}C^k_n\varphi_c^{2n-2k}(\partial_x\varphi_c)^{2k+1}\Big)\phi\Big]dx
dt \\
&=-\int_0^{+\infty}\int_{\R}\phi\cdot\partial_x p\ast\Big[\sum_{k=1}^n\frac{(-1)^{k-1}}{2k-1}C^k_n\varphi_c^{2n-2k+1}(\partial_x\varphi_c)^{2k}\\
&\qquad\qquad\qquad\qquad\qquad\qquad+\frac{1}{2n+1}\Big(2n+\sum_{k=1}^n\frac{(-1)^{k+1}}{2k+1}C^k_n\Big)\varphi_c^{2n+1}\Big]dxdt.
\end{aligned}
\end{eqnarray}
Note that $\partial_xp(x)=-\mbox{sign}(x)e^{-|x|}/2\ \mbox{for} \ x\in \R$, we deduce
\begin{eqnarray}\label{2.9}
\begin{aligned}
&\partial_x p\ast\Big[\sum_{k=1}^n\frac{(-1)^{k-1}}{2k-1}C^k_n\varphi_c^{2n-2k+1}(\partial_x\varphi_c)^{2k}+\frac{1}{2n+1}
\Big(2n+\sum_{k=1}^n\frac{(-1)^{k+1}}{2k+1}C^k_n\Big)\varphi_c^{2n+1}\Big]\\
&=-\frac{1}{2}\int^{+\infty}_{-\infty}\mbox{sign}(x-y)e^{-|x-y|}\Big[\frac{1}{2n+1}
\Big(2n+\sum_{k=1}^n\frac{(-1)^{k+1}}{2k+1}C^k_n\Big)\\
&\qquad\qquad\qquad\qquad\qquad\qquad\qquad+\sum_{k=1}^n\frac{(-1)^{k-1}}{2k-1}C^k_n\mbox{sign}^{2k}(y-ct)
 \Big]a^{2n+1}e^{-(2n+1)|y-ct|}dy.
\end{aligned}
\end{eqnarray}
For $x>ct$, we can split the right hand side of (\ref{2.9}) into
the following three parts,
\begin{eqnarray}\label{2.10}
&&\partial_x p\ast\Big[\sum_{k=1}^n\frac{(-1)^{k-1}}{2k-1}C^k_n\varphi_c^{2n-2k+1}(\partial_x\varphi_c)^{2k}+\frac{1}{2n+1}
\Big(2n+\sum_{k=1}^n\frac{(-1)^{k+1}}{2k+1}C^k_n\Big)\varphi_c^{2n+1}\Big]\nonumber\\
&=&-\frac{1}{2}\Big[\frac{1}{2n+1}
\Big(2n+\sum_{k=1}^n\frac{(-1)^{k+1}}{2k+1}C^k_n\Big)+\sum_{k=1}^n\frac{(-1)^{k-1}}{2k-1}C^k_n\Big]a^{2n+1}\nonumber\\
&&\times
\Big[\int^{ct}_{-\infty}+\int_{ct}^x+\int_x^{+\infty}\Big]\mbox{sign}(x-y)e^{-|x-y|}e^{-(2n+1)|y-ct|}dy\nonumber\\
 &\triangleq&I_1+I_2+I_3.
\end{eqnarray}
A direct calculation for each one of the terms $I_i,1\leq i \leq 3,$ yields
\begin{eqnarray*}
I_1&=&-\frac{1}{2}\Big[\frac{1}{2n+1}
\Big(2n+\sum_{k=1}^n\frac{(-1)^{k+1}}{2k+1}C^k_n\Big)+\sum_{k=1}^n\frac{(-1)^{k-1}}{2k-1}C^k_n\Big]a^{2n+1}\times
\int^{ct}_{-\infty}e^{-(x-y)}e^{(2n+1)(y-ct)}dy\nonumber\\
&=&-\frac{1}{2}\Big[\frac{1}{2n+1}
\Big(2n+\sum_{k=1}^n\frac{(-1)^{k+1}}{2k+1}C^k_n\Big)+\sum_{k=1}^n\frac{(-1)^{k-1}}{2k-1}C^k_n\Big]a^{2n+1}\times e^{-[x+(2n+1)ct]}
\int^{ct}_{-\infty}e^{(2n+2)y}dy\nonumber\\
&=&-\frac{1}{4n+4}\Big[\frac{1}{2n+1}
\Big(2n+\sum_{k=1}^n\frac{(-1)^{k+1}}{2k+1}C^k_n\Big)+\sum_{k=1}^n\frac{(-1)^{k-1}}{2k-1}C^k_n\Big]a^{2n+1}e^{ct-x},
\end{eqnarray*}
\begin{eqnarray*}
I_2&=&-\frac{1}{2}\Big[\frac{1}{2n+1}
\Big(2n+\sum_{k=1}^n\frac{(-1)^{k+1}}{2k+1}C^k_n\Big)+\sum_{k=1}^n\frac{(-1)^{k-1}}{2k-1}C^k_n\Big]a^{2n+1}\nonumber\\
&&\times
\int_{ct}^xe^{-(x-y)}e^{-(2n+1)(y-ct)}dy\nonumber\\
&=&-\frac{1}{2}\Big[\frac{1}{2n+1}
\Big(2n+\sum_{k=1}^n\frac{(-1)^{k+1}}{2k+1}C^k_n\Big)+\sum_{k=1}^n\frac{(-1)^{k-1}}{2k-1}C^k_n\Big]a^{2n+1}\times e^{-[x-(2n+1)ct]}
\int_{ct}^x e^{-2ny}dy\nonumber\\
&=&-\frac{1}{4n}\Big[\frac{1}{2n+1}
\Big(2n+\sum_{k=1}^n\frac{(-1)^{k+1}}{2k+1}C^k_n\Big)+\sum_{k=1}^n\frac{(-1)^{k-1}}{2k-1}C^k_n\Big]\times a^{2n+1}\big(e^{ct-x}-e^{(2n+1)(ct-x)}\big),
\end{eqnarray*}
and
\begin{eqnarray*}
I_3&=&\frac{1}{2}\Big[\frac{1}{2n+1}
\Big(2n+\sum_{k=1}^n\frac{(-1)^{k+1}}{2k+1}C^k_n\Big)+\sum_{k=1}^n\frac{(-1)^{k-1}}{2k-1}C^k_n\Big]a^{2n+1}\times
\int_x^{+\infty}e^{x-y}e^{-(2n+1)(y-ct)}dy\nonumber\\
&=&\frac{1}{2}\Big[\frac{1}{2n+1}
\Big(2n+\sum_{k=1}^n\frac{(-1)^{k+1}}{2k+1}C^k_n\Big)+\sum_{k=1}^n\frac{(-1)^{k-1}}{2k-1}C^k_n\Big]a^{2n+1}\times e^{[x+(2n+1)ct]}
\int_x^{+\infty} e^{-(2n+2)y}dy\nonumber\\
&=&\frac{1}{4n+4}\Big[\frac{1}{2n+1}
\Big(2n+\sum_{k=1}^n\frac{(-1)^{k+1}}{2k+1}C^k_n\Big)+\sum_{k=1}^n\frac{(-1)^{k-1}}{2k-1}C^k_n\Big]a^{2n+1}e^{(2n+1)(ct-x)}.
\end{eqnarray*}
Plugging the above equalities $I_1$-$I_3$ into (\ref{2.10}), we find that for $x>ct$,
\begin{eqnarray}\label{2.11}
&&\partial_x p\ast\Big[\sum_{k=1}^n\frac{(-1)^{k-1}}{2k-1}C^k_n\varphi_c^{2n-2k+1}(\partial_x\varphi_c)^{2k}+\frac{1}{2n+1}
\Big(2n+\sum_{k=1}^n\frac{(-1)^{k+1}}{2k+1}C^k_n\Big)\varphi_c^{2n+1}\Big]\nonumber\\
&=&-\frac{1}{4n\cdot(n+1)}\Big[
\Big(2n+\sum_{k=1}^n\frac{(-1)^{k+1}}{2k+1}C^k_n\Big)+(2n+1)\sum_{k=1}^n\frac{(-1)^{k-1}}{2k-1}C^k_n\Big]\nonumber\\
&&\times a^{2n+1}\big(e^{ct-x}-e^{(2n+1)(ct-x)}\big).
\end{eqnarray}
While for $x\leq ct$, we also split the right hand side of (\ref{2.9}) into three parts,
\begin{eqnarray}\label{2.12}
&&\partial_x p\ast\Big[\sum_{k=1}^n\frac{(-1)^{k-1}}{2k-1}C^k_n\varphi_c^{2n-2k+1}(\partial_x\varphi_c)^{2k}+\frac{1}{2n+1}
\Big(2n+\sum_{k=1}^n\frac{(-1)^{k+1}}{2k+1}C^k_n\Big)\varphi_c^{2n+1}\Big]\nonumber\\
&=&-\frac{1}{2}\Big[\frac{1}{2n+1}
\Big(2n+\sum_{k=1}^n\frac{(-1)^{k+1}}{2k+1}C^k_n\Big)+\sum_{k=1}^n\frac{(-1)^{k-1}}{2k-1}C^k_n\Big]a^{2n+1}\nonumber\\
&&\times
\Big[\int^{x}_{-\infty}+\int_{x}^{ct}+\int_{ct}^{+\infty}\Big]\mbox{sign}(x-y)e^{-|x-y|}e^{-(2n+1)|y-ct|}dy\nonumber\\
&\triangleq&II_1+II_2+II_3.
\end{eqnarray}
Applying a similar computation as $I_i,1\leq i \leq 3,$ to the terms
$II_1$-$II_3$ on the right hand side of (\ref{2.12}), we find that for $x\leq ct$,
\begin{eqnarray}\label{2.13}
&&\partial_x p\ast\Big[\sum_{k=1}^n\frac{(-1)^{k-1}}{2k-1}C^k_n\varphi_c^{2n-2k+1}(\partial_x\varphi_c)^{2k}+\frac{1}{2n+1}
\Big(2n+\sum_{k=1}^n\frac{(-1)^{k+1}}{2k+1}C^k_n\Big)\varphi_c^{2n+1}\Big]\nonumber\\
&=&\frac{1}{4n\cdot(n+1)}\Big[
\Big(2n+\sum_{k=1}^n\frac{(-1)^{k+1}}{2k+1}C^k_n\Big)+(2n+1)\sum_{k=1}^n\frac{(-1)^{k-1}}{2k-1}C^k_n\Big]\nonumber\\
&&\times a^{2n+1}\big(e^{x-ct}-e^{(2n+1)(x-ct)}\big).
\end{eqnarray}
Now we claim that the following identity holds,
\begin{eqnarray}\label{2.14}
\frac{1}{4n\cdot(n+1)}\Big[
\Big(2n+\sum_{k=1}^n\frac{(-1)^{k+1}}{2k+1}C^k_n\Big)+(2n+1)\sum_{k=1}^n\frac{(-1)^{k-1}}{2k-1}C^k_n\Big]=1-\sum_{k=1}^n\frac{(-1)^{k+1}}{2k+1}C^k_n.
\end{eqnarray}
Indeed, to prove (\ref{2.14}), we need to verify the following identity:
\begin{eqnarray}\label{2.15}
(2n+1)\sum_{k=1}^{n-1}\frac{(-1)^{k+1}}{2k+1}C^k_n+\sum_{k=1}^n\frac{(-1)^{k-1}}{2k-1}C^k_n=2n+(-1)^n.
\end{eqnarray}
Indeed, using the combination formula: $\sum\limits_{k=0}^{n}\frac{(-1)^{k}}{2k+1}C^k_n=\frac{(2n)!!}{(2n+1)!!}$, we have
\begin{eqnarray}\label{2.16}
\sum_{k=1}^{n-1}\frac{(-1)^{k+1}}{2k+1}C^k_n=-\Big(\frac{(2n)!!}{(2n+1)!!}-1-\frac{(-1)^n}{2n+1}\Big).
\end{eqnarray}
In contrast with (\ref{2.16}) and (\ref{2.15}), it remains to show
\begin{eqnarray}\label{2.17}
\sum_{k=1}^n\frac{(-1)^{k-1}}{2k-1}C^k_n=\frac{(2n)!!}{(2n-1)!!}-1.
\end{eqnarray}
We prove it by the induction method. Clearly,  (\ref{2.17}) holds for $n=1$. By induction, assuming (\ref{2.17}) holds for $n$, we then have
\begin{eqnarray*}
&&\sum_{k=1}^{n+1}\frac{(-1)^{k-1}}{2k-1}C^k_{n+1}=\sum_{k=1}^{n}\frac{(-1)^{k-1}}{2k-1}\big(C^k_{n}+C^{k-1}_{n}\big)+\frac{(-1)^n}{2n+1}\\
&=&\sum_{k=1}^{n}\frac{(-1)^{k-1}}{2k-1}C^k_{n}+\sum_{k=1}^{n}\frac{(-1)^{k-1}}{2k-1}C^{k-1}_{n}+\frac{(-1)^n}{2n+1}\\
&=&\frac{(2n)!!}{(2n-1)!!}-1+\sum_{k=0}^{n}\frac{(-1)^{k}}{2k+1}C^{k}_{n}=\frac{(2n+2)!!}{(2n+1)!!}-1,
\end{eqnarray*}
which implies (\ref{2.17}). Hence, by (\ref{T1}), (\ref{T2}), (\ref{2.11}), (\ref{2.13}) and (\ref{2.14}), we deduce that for all $(t,x)\in (0,+\infty)\times \R,$
\begin{eqnarray}\label{2.18}
&&\mbox{sign}(x-ct)\varphi_c
\Big[c-\Big(1-\sum_{k=1}^n\frac{(-1)^{k+1}}{2k+1}C^k_n\Big)\varphi_c^{2n}\Big]+\partial_x p\ast\Big[\sum_{k=1}^n\frac{(-1)^{k-1}}{2k-1}C^k_n\varphi_c^{2n-2k+1}(\partial_x\varphi_c)^{2k}\nonumber\\
&&+\frac{1}{2n+1}
\Big(2n+\sum_{k=1}^n\frac{(-1)^{k+1}}{2k+1}C^k_n\Big)\varphi_c^{2n+1}\Big](t,x)=0.
\end{eqnarray}
Thanks to (\ref{2.6}), (\ref{2.8}) and (\ref{2.18}), we conclude that
\begin{eqnarray*}
\begin{aligned}
&\int_0^{+\infty}\int_{\R}\Big[\varphi_c\partial_t\phi+\frac{1}{2n+1}\varphi_c^{2n+1}\partial_x\phi
+\Big(\sum_{k=1}^n\frac{(-1)^{k+1}}{2k+1}C^k_n\varphi_c^{2n-2k}(\partial_x\varphi_c)^{2k+1}\Big)\phi\\
&\qquad\qquad\quad+(1-\partial_x^2)^{-1}\Big(\frac{2n}{2n+1}\varphi_c^{2n+1}+
\sum_{k=1}^n\frac{(-1)^{k-1}}{2k-1}C^k_n\varphi_c^{2n-2k+1}(\partial_x\varphi_c)^{2k}\Big)\partial_x\phi\\
&\qquad\qquad\qquad-(1-\partial_x^2)^{-1}\Big(\sum_{k=1}^n\frac{(-1)^{k+1}}{2k+1}C^k_n\varphi_c^{2n-2k}(\partial_x\varphi_c)^{2k+1}\Big)\phi
\Big]dxdt\\
&\qquad\qquad\qquad\qquad+\int_{\R}\varphi_{0, c}(x)\phi(0,x)dx=0,
\end{aligned}
\end{eqnarray*}
for any $\phi(t,x)\in C_c^{\infty}([0,\infty)\times \R)$. This completes the proof of Theorem 2.1.
\end{proof}

It was shown in \cite{R-A} that the gmCH equation \eqref{1.1} can be written as the Hamiltonian form
\begin{eqnarray}\label{gmCH-Ham}
y_t=-D_x \frac{\delta H_{gmCH}}{\delta u},
\end{eqnarray}
where
\begin{eqnarray*}
H_{gmCH}=\frac 1 {2(n+1)}\int_{\mathbb R}u(u^2-u_x^2)^n y dx.
\end{eqnarray*}
It is easy to check that $H_{gmCH}$ is equivalent to
\begin{eqnarray}\label{gmCH-Ham-1}
F(u)= \int_{\R}\left(u^{2n+2}+\sum_{k=1}^{n}\frac{(-1)^{k+1}}{2k-1}C^k_{n+1}u^{2n-2k+2}u_x^{2k}+(-1)^n\frac{u^{2n+2}_x}{2n+1}\right)dx.
\end{eqnarray}
Actually, we can also directly prove that $F(u)$ is a conservation law of equation \eqref{1.1} (See Appendix A.3.). As for the CH and modified CH equation, the gmCH equation also admits the $H^1$ conserved density
\begin{eqnarray}\label{gmCH-Ham-2}
E(u)= \int_{\R}\big(u^2+u^2_x\big)dx.
\end{eqnarray}

\section{Stability of peakons}

\newtheorem{theorem3}{Theorem}[section]
\newtheorem{lemma3}{Lemma}[section]
\newtheorem {remark3}{Remark}[section]
\newtheorem {definition3}{Definition}[section]
\newtheorem{corollary3}{Corollary}[section]
\par

In this section, we prove  orbital stability of singe peakons for the gmCH equation (\ref{1.1}).
We now present the stability result.
\begin{theorem3}
The peakon $\varphi_c(t,x)$ defined in (\ref{2.2}) traveling
with the speed $c>0$ is orbitally stable in the following sense. If $u_0(x)\in H^s(\R)$, for some $s>5/2,$
$y_0(x)=(1-\partial^2_x)u_0\not\equiv 0$ is nonnegative, and
\begin{eqnarray*}
\big\|u(0,\cdot)-\varphi_c\big\|_{H^{1}(\R)}<\varepsilon, \quad\quad \mbox{for} \quad 0<\varepsilon<\big(3-2\sqrt{2}\big)a,
\end{eqnarray*}
where $a>0$ is defined by $c=\left(1-\sum_{k=1}^n\frac{(-1)^{k+1}}{2k+1}C^k_n\right)a^{2n}$ according to Theorem 2.1.
Then the corresponding solution $u(t,x)$ of equation (\ref{1.1}) satisfies
\begin{eqnarray*}
\sup_{t\in [0,T)}\big\|u(t,\cdot)-\varphi_c(\cdot-\xi(t))\big\|_{H^1(\R)}\lesssim \sqrt{3a\varepsilon+4a\sqrt{A(c,\|u_0\|_{H^s})\varepsilon}},
\end{eqnarray*}
where $T>0$ is the maximal existence time, $\xi(t) \in {\R}$ is the maximum point of function
$u(t,\cdot)$ and the constant $A(n,c,\|u_0\|_{H^s})>0$ depends only on wave speed $c>0$, integer $n>0$ and the norm $\|u_0\|_{H^s}$.
\end{theorem3}
\begin{remark2}
Note that by the above stability theorem we mean that even if a solution $u$ which is initially close to a peakon $\varphi_c$ blows up
in a finite time, it stays close to some translate of the peakon up to the breaking time.
\end{remark2}
We break the proof of Theorem 3.1 into several lemmas. Note that the assumptions on the initial
profile in Theorem 3.1 guarantee the existence of the unique local positive solution of equation (\ref{1.1}) by Lemmas 2.1-2.2. It is obvious that
$\varphi_c(x)=a\varphi(x)=ae^{-|x|}\in H^1(\R)$ has the peak at $x=0,$
and hence
\begin{eqnarray*}
\max\limits_{x\in\R}\big\{\varphi_c(x)\big\}=\varphi_c(0)=a.
\end{eqnarray*}
By a simple computation, we have
\begin{eqnarray}\label{3.2}
E(\varphi_c)=\|\varphi_c\|_{H^1}^2=a^2\int_{\R}(\varphi^2+\varphi^2_x)dx=2a^2,
\end{eqnarray}
and
\begin{eqnarray*}
F(\varphi_c)&=&a^{2n+2}\int_{\R}\left(\varphi^{2n+2}+\sum_{k=1}^{n}\frac{(-1)^{k+1}}{2k-1}C^k_{n+1}\varphi^{2n-2k+2}\varphi_x^{2k}
+(-1)^n\frac{\varphi^{2n+2}_x}{2n+1}\right)dx\nonumber\\
&=&\frac{a^{2n+2}}{n+1} \left(1+\sum_{k=1}^{n+1}\frac{(-1)^{k+1}}{2k-1}C^k_{n+1}\right),
\end{eqnarray*}
where $a$ is given implicitly by $c=\left(1-\sum_{k=1}^n\frac{(-1)^{k+1}}{2k+1}C^k_n\right)a^{2n}$ according to Theorem 2.1.

\begin{lemma3}
For every $u\in H^1(\R)$ and $\xi\in\R,$ we have
\begin{eqnarray}\label{3.4}
E(u)-E\big(\varphi_c(\cdot-\xi)\big)=\big\|u-\varphi_c(\cdot-\xi)\big\|^2_{H^1}
+4 a\big(u(\xi)-a\big),
\end{eqnarray}
where $a$ is defined by $c=\left(1-\sum_{k=1}^{n}\frac{(-1)^{k+1}}{2k+1}C^k_n\right)a^{2n}$.
\end{lemma3}
\begin{proof}
Note the relation $\varphi-\partial_x^2\varphi=2\delta$. Here $\delta$ denotes the Dirac distribution.
For simplicity, we abuse notation by writing integrals instead of the $H^{-1}/H^1$ duality pairing.
Hence we have
\begin{eqnarray*}
\big\|u-\varphi_c(\cdot-\xi)\big\|^2_{H^1}
&=&\int_{\R}(u^2+u^2_x)dx+\int_{\R}\big(\varphi_c^2+(\partial_x\varphi_c)^2\big)dx
-2a\int_{\R}u_x(x)\varphi_x(x-\xi)dx\\
&&-2a\int_{\R}u(x)\varphi (x-\xi)dx\nonumber\\
&=&E(u)+E\big(\varphi_c(\cdot-\xi)\big)-2a\int_{\R}(1-\partial^2_x)\varphi (x-\xi)u(x)dx\nonumber\\
&=&E(u)+E\big(\varphi_c(\cdot-\xi)\big)-4a\int_{\R}\delta (x-\xi)u(x)dx\nonumber\\
&=&E(u)+E\big(\varphi_c(\cdot-\xi)\big)-4au(\xi)=E(u)-E\big(\varphi_c(\cdot-\xi)\big)-4a\big(u(\xi)-a\big),
\end{eqnarray*}
where we used integration by parts and (\ref{3.2}). This completes the proof of Lemma 3.1.
\end{proof}

Next, we derive a crucial polynomial inequality relating the two conserved
quantities $E(u)$ and $F(u)$ to the maximal value of approximate solutions.
\begin{lemma3}
Assume $u_0\in H^s(\R),s>5/2$, and $y_0\geq 0.$ Let $u(t, x)$ be the positive solution of the Cauchy problem of the gmCH equation (\ref{1.1}) with
initial data $u_0$. Denote $M(t)\triangleq\max_{x\in\R}\{u(t, x)\}$. Then
\begin{eqnarray}\label{3.5}
\frac{n(2-c_1)}{n+1}M^{2n+2}(t)-\frac{2-c_1}{2}M^{2n}(t)E(u)+F(u)\leq 0.
\end{eqnarray}
\end{lemma3}

\begin{proof}
Let $M(t)$ be taken at $x=\xi(t)$, and define the same function $g$ as in \cite{Constantin-S}
\begin{equation*}
g(t, x)\triangleq\left\{\begin{array}{ll}u(t, x)- u_x(t, x), & x<\xi(t),\\
u(t, x)+ u_x(t, x),& x>\xi(t).
\end{array}\right.
\end{equation*}
Then we have
\begin{eqnarray}\label{3.6}
\int_{\R}g^2(t, x)dx=E(u)-2M^2(t).
\end{eqnarray}
Next we introduce the following function
\begin{eqnarray*}
h(t, x)\triangleq\left\{\begin{array}{ll}
\left(u^{2n}+\sum\limits_{k=1}^{2n-1}c_ku^{2n-k}u_x^k+(-1)^n\frac{1}{2n+1}u^{2n}_x\right)(t, x), & x<\xi(t),\\
\left(u^{2n}+\sum\limits_{k=1}^{2n-1}d_ku^{2n-k}u_x^k+(-1)^n\frac{1}{2n+1}u^{2n}_x\right)(t, x),& x>\xi(t).
\end{array}\right.
\end{eqnarray*}
where $c_k,d_k,k=1,2,\cdot\cdot\cdot,2n-1$, are constants given by
\begin{equation}\label{3.7}
\left\{\begin{array}{lll}c_1=-d_1=\frac{1}{2}+\sum\limits_{j=1}^{n+1}(-1)^{j+1}\frac{2j-3}{2(2j-1)}C^j_{n+1},\\
c_{2m}=d_{2m}=\sum\limits_{j=m+1}^{n+1}(-1)^{j+1}\frac{2j-(2m+1)}{2j-1}C^j_{n+1},& m=1,2,\cdot\cdot\cdot,n-1,\\
c_{2m-1}=-d_{2m-1}=\sum\limits_{j=m+1}^{n+1}(-1)^{j+1}\frac{2j-2m}{2j-1}C^j_{n+1},& m=2,3\cdot\cdot\cdot,n.
\end{array}\right.
\end{equation}
Integrating by parts, we calculate
\begin{eqnarray}\label{3.8}
&&\int_{\R}h(t, x)g^2(t, x)dx\nonumber\\
&=&\int_{-\infty}^{\xi}\big(u- u_x\big)^2\big(u^{2n}+\sum\limits_{k=1}^{2n-1}c_ku^{2n-k}u_x^k+(-1)^n\frac{1}{2n+1}u^{2n}_x
\big)dx\nonumber\\
&&+\int^{\infty}_{\xi}\big(u+u_x\big)^2\big(u^{2n}+\sum\limits_{k=1}^{2n-1}d_ku^{2n-k}u_x^k+(-1)^n\frac{1}{2n+1}u^{2n}_x\big)dx\nonumber\\
&=&\int_{-\infty}^{\xi}\big(u^{2n+2}+\sum_{k=1}^{n}\frac{(-1)^{k+1}}{2k-1}C^k_{n+1}u^{2n-2k+2}u_x^{2k}+(-1)^n\frac{u^{2n+2}_x}{2n+1}\big)dx\nonumber\\
&&+(c_1-2)\int_{-\infty}^{\xi}u^{2n+1}u_xdx +\int^{\infty}_{\xi}\big(u^{2n+2}+\sum_{k=1}^{n}\frac{(-1)^{k+1}}{2k-1}C^k_{n+1}u^{2n-2k+2}u_x^{2k}\nonumber\\
&&+(-1)^n\frac{u^{2n+2}_x}{2n+1}\big)dx+(d_1+2)\int^{\infty}_{\xi}u^{2n+1}u_xdx\nonumber\\
&=&F(u)+\frac{c_1-2}{n+1}u^{2n+2}(t, \xi)=F(u)-\frac{2-c_1}{n+1}M^{2n+2}(t),
\end{eqnarray}
where the last third equality in (\ref{3.8}) will be shown in Appendix A.1.\\
On the other hand, by Lemma 2.2, the solution $u$ satisfies
\begin{eqnarray}\label{3.9}
u(t,x)\geq 0,\ \mbox{and} \ (u\pm u_x)(t,x)\geq 0 \quad \mbox{for}\ \forall (t,x)\in [0,T)\times \R.
\end{eqnarray}
We now claim that
\begin{eqnarray}\label{3.10}
h(t,x)\leq \frac{2-c_1}{2}u^{2n}(t,x) \quad  \mbox{for}\ \forall (t,x)\in [0,T)\times \R.
\end{eqnarray}
In view of the expression for $h$, it suffices to show that
\begin{eqnarray}\label{3.11}
\sum\limits_{k=1}^{2n-1}c_ku^{2n-k}u_x^k\leq -\frac{c_1}{2}u^{2n}-\frac{(-1)^n}{2n+1}u^{2n}_x,
\end{eqnarray}
and
\begin{eqnarray}\label{3.12}
\sum\limits_{k=1}^{2n-1}d_ku^{2n-k}u_x^k\leq -\frac{c_1}{2}u^{2n}-\frac{(-1)^n}{2n+1}u^{2n}_x.
\end{eqnarray}
By (\ref{3.9}), and denoting $z=\frac{u_x}{u}$,
we find that the proof of (\ref{3.11}) is equivalent to prove the non-positivity of the following function $f(z)$. That is,
\begin{eqnarray}\label{0.1}
f(z)=\frac{c_{2n-1}}{2}z^{2n}+\sum\limits_{k=1}^{2n-1}c_kz^k+\frac{c_1}{2}\leq 0, \quad \mbox{for} \ z \in [-1,1].
\end{eqnarray}
By the fact $c_{2j+1}-2c_{2j}+c_{2j-1}=0\ (j=1,2,\cdot\cdot\cdot,n-1)$ in (\ref{4.3}), we can rewrite the above function
$f(z)$ as follows.
\begin{eqnarray*}
f(z)&=&\frac{c_{2n-1}}{2}z^{2n}+c_{2n-1}z^{2n-1}+\sum\limits_{j=1}^{n-1}\big(c_{2j}z^{2j}+c_{2j-1}z^{2j-1}
\big)+\frac{c_1}{2}\nonumber\\
&=& \frac{c_{2n-1}}{2}z^{2n}+c_{2n-1}z^{2n-1}+\sum\limits_{j=1}^{n-1}\big(\frac{c_{2j+1}}{2}z^{2j}+
\frac{c_{2j-1}}{2}z^{2j}+c_{2j-1}z^{2j-1}
\big)+\frac{c_1}{2}\nonumber\\
&=& \sum\limits_{k=1}^{n}\big(\frac{c_{2k-1}}{2}z^{2k}+
c_{2k-1}z^{2k-1}+\frac{c_{2k-1}}{2}z^{2k-2}\big)=\frac{(1+z)^2}{2}\sum\limits_{k=1}^{n}c_{2k-1}z^{2k-2}.
\end{eqnarray*}
Hence, to prove (\ref{0.1}), it suffices to show that
\begin{eqnarray}\label{0}
\phi(z)=\sum\limits_{k=1}^{n} c_{2k-1}z^{2k-2}\leq 0, \quad \mbox{for} \ z \in [-1,1].
\end{eqnarray}
Note that $\phi(z)$ is an even function and continuous at $z=1$
, we only need to prove
the above inequality (\ref{0}) holds for $z\in [0,1).$
Using the formula $c_{2k-1}\ (k=1,2,\cdot\cdot\cdot,n)$ in (\ref{3.7}), and exchanging the order of summation,
we obtain
\begin{eqnarray}\label{0.2}
\phi(z)&=&c_1+\sum\limits_{k=2}^{n} \sum\limits_{j=k+1}^{n+1}(-1)^{j+1}\frac{2j-2k}{2j-1}C_{n+1}^jz^{2k-2}\nonumber\\
&=& c_1+\sum\limits_{j=3}^{n+1} \sum\limits_{k=2}^{j-1}(-1)^{j+1}\big(1+\frac{1-2k}{2j-1}\big)C_{n+1}^jz^{2k-2}\nonumber\\
&\triangleq& c_1+\tilde{\phi}_1(z)+\tilde{\phi}_2(z).
\end{eqnarray}
For $\tilde{\phi}_1(z)$, we have
\begin{eqnarray}\label{0.3}
\tilde{\phi}_1(z)&=&\sum\limits_{j=3}^{n+1} (-1)^{j+1}C_{n+1}^j\big(\frac{z^{2}(1-z^{2(j-2)})}{1-z^2}\big)
\nonumber\\
&=&\frac{z^2}{1-z^2}\sum\limits_{j=3}^{n+1} (-1)^{j+1}C_{n+1}^j+\frac{z^{-2}}{1-z^2}\sum\limits_{j=3}^{n+1} (-1)^{j}C_{n+1}^jz^{2j}\nonumber\\
&\triangleq& \tilde{\phi}_{1,1}(z)+\tilde{\phi}_{1,2}(z).
\end{eqnarray}
A direct calculation for each one of the terms $\tilde{\phi}_{1,i}(z),$ $i=1,2,$ yields
\begin{eqnarray}\label{0.4}
\tilde{\phi}_{1,1}(z)&=&-\frac{z^2}{1-z^2}\big(\sum\limits_{j=0}^{n+1} (-1)^{j}C_{n+1}^j
-\sum\limits_{j=0}^{2} (-1)^{j}C_{n+1}^j\big)\nonumber\\
&=&\frac{z^2}{1-z^2}\cdot\sum\limits_{j=0}^{2} (-1)^{j}C_{n+1}^j=\frac{n^2-n}{2} \frac{z^2}{1-z^2},
\end{eqnarray}
and
\begin{eqnarray}\label{0.5}
\tilde{\phi}_{1,2}(z)&=&\frac{z^{-2}}{1-z^2}
\big(\sum\limits_{j=0}^{n+1} (-z^2)^{j}C_{n+1}^j
-\sum\limits_{j=0}^{2}(-z^2)^{j}C_{n+1}^j\big)\nonumber\\
&=&\frac{z^{-2}}{1-z^2}
\big((1-z^2)^{n+1}-1+(n+1)z^2-z^4C_{n+1}^2\big) \triangleq \frac{z^{-2}}{1-z^2}\cdot\omega(z).
\end{eqnarray}
Thus, substituting (\ref{0.4})-(\ref{0.5}) into (\ref{0.3}), we get
\begin{eqnarray}\label{0.6}
\tilde{\phi}_1(z)=\frac{n^2-n}{2}\frac{z^2}{1-z^2}+\frac{\omega(z)z^{-2}}{1-z^2}.
\end{eqnarray}
For $\tilde{\phi}_2(z)$, we have
\begin{eqnarray}\label{0.7}
\tilde{\phi}_{2}(z)&=&\frac{d}{dz}\Big(\sum\limits_{j=3}^{n+1} \sum\limits_{k=2}^{j-1}
\frac{(-1)^{j}C_{n+1}^j}{2j-1}z^{2k-1}\Big)\nonumber\\
&=&\frac{d}{dz}\Big(\sum\limits_{j=3}^{n+1} \frac{(-1)^{j}C_{n+1}^j}{2j-1}\cdot\frac{z^3(1-z^{2(j-2)})}{1-z^2}\Big)\nonumber\\
&=& \frac{d}{dz}\Big( \frac{z^3}{1-z^2}\cdot\sum\limits_{j=3}^{n+1} \frac{(-1)^{j}C_{n+1}^j}{2j-1}
-\frac{1}{1-z^2}\cdot\sum\limits_{j=3}^{n+1} \frac{(-1)^{j}C_{n+1}^j}{2j-1} z^{2j-1}\Big)\nonumber\\
&\triangleq& \frac{d}{dz}\big(\tilde{\Phi}_{2,1}(z)+\tilde{\Phi}_{2,2}(z)\big).
\end{eqnarray}
Applying a similar computation as in (\ref{0.5}) to $\tilde{\Phi}_{2,i}(z),i=1,2,$ we obtain
\begin{eqnarray}\label{0.8}
\tilde{\Phi}_{2,1}(z)&=& \frac{z^3}{1-z^2}\cdot \int_0^1\sum\limits_{j=3}^{n+1} (-1)^{j}C_{n+1}^j z^{2j-2}dz
\nonumber\\
&=&\frac{z^3}{1-z^2}\cdot \int_0^1\sum\limits_{j=3}^{n+1}C_{n+1}^j (-z^2)^{j}\frac{dz}{z^2}=\frac{z^3}{1-z^2}\cdot \int_0^1 \frac{\omega(z)}{z^2}dz,
\end{eqnarray}
and
\begin{eqnarray}\label{0.9}
\tilde{\Phi}_{2,2}(z)&=& \frac{-1}{1-z^2}\cdot \int_0^z\sum\limits_{j=3}^{n+1} (-z^2)^{j}C_{n+1}^j \frac{dz}{z^2}=\frac{-1}{1-z^2}\cdot \int_0^z \frac{\omega(z)}{z^2}dz
\end{eqnarray}
Thus, substituting (\ref{0.8})-(\ref{0.9}) into (\ref{0.7}), we get
\begin{eqnarray}\label{0.10}
\tilde{\phi}_2(z)=
\frac{3z^2\cdot\int_0^1 \frac{\omega(z)}{z^2}dz-\omega(z)z^{-2}}{1-z^2}+
\frac{2z(z^3\cdot \int_0^1 \frac{\omega(z)}{z^2}dz-\int_0^z \frac{\omega(z)}{z^2}dz)}{(1-z^2)^2}.
\end{eqnarray}
For $c_1$, similarly as in (\ref{0.8}), by (\ref{3.7}), we have
\begin{eqnarray}\label{0.11}
c_1&=&\frac{1}{2}+\sum\limits_{j=1}^{n+1}(-1)^{j+1}\frac{C^j_{n+1}}{2}-\sum\limits_{j=1}^{n+1}(-1)^{j+1}\frac{C^j_{n+1}}{2j-1}\nonumber\\
&=&1+\sum\limits_{j=0}^{n+1}(-1)^{j+1}\frac{C^j_{n+1}}{2}+\sum\limits_{j=1}^{n+1}\frac{(-1)^{j}C^j_{n+1}}{2j-1}=-n+\frac{C_{n+1}^2}{3}+ \int_0^1 \frac{\omega(z)}{z^2}dz.
\end{eqnarray}
Note that,
\begin{eqnarray}\label{0.12}
\int_0^1 \frac{\omega(z)}{z^2}dz&=&-\int_0^1 \sum\limits_{k=0}^{n}(1-z^2)^k dz+\int_0^1 \big(n+1-z^2C^2_{n+1}\big)dz\nonumber\\
&=&-\sum\limits_{k=1}^{n}\frac{(2k)!!}{(2k+1)!!}+n-\frac{C^2_{n+1}}{3}\triangleq -B+n-\frac{C^2_{n+1}}{3}.
\end{eqnarray}
Thus, combining (\ref{0.2}), (\ref{0.6}) with (\ref{0.10})-(\ref{0.12}), we obtain
\begin{eqnarray*}
\phi(z)&=&-B+\frac{(2n-3B)z^2}{1-z^2}+\frac{2(-B+n-\frac{C^2_{n+1}}{3})z^4}{(1-z^2)^2}\nonumber\\
&&-\frac{2z\cdot\int_0^z(n+1-z^2C^2_{n+1})dz}{(1-z^2)^2}+\frac{2z\cdot\int_0^z \sum_{k=0}^{n}(1-z^2)^k dz}{(1-z^2)^2}\nonumber\\
&=&\frac{-B-Bz^2-2z^2}{(1-z^2)^2}+\frac{2z^2+2z \int_0^z \sum_{k=1}^{n}(1-z^2)^k dz}{(1-z^2)^2}\nonumber\\
&\leq&\frac{-B-Bz^2}{(1-z^2)^2}+\frac{2z\cdot\int_0^1 \sum_{k=1}^{n}(1-z^2)^k dz}{(1-z^2)^2}=-\frac{B}{(1+z)^2}\leq 0,
\end{eqnarray*}
which gives the desired result (\ref{0}) for $z\in [0,1)$. Then by (\ref{3.7}) and (\ref{4.4}), using a similar argument, we find that (\ref{3.12}) is equivalent to $\sum_{k=1}^{n} c_{2k-1}z^{2k-2}(z-1)^2\leq 0,$ which is obviously true by (\ref{0}). Hence, we complete the proof of the claim (\ref{3.10}). Therefore, in view of (\ref{3.6}), (\ref{3.8}) and (\ref{3.10}), we conclude that
\begin{eqnarray*}
&&F(u)-\frac{2-c_1}{n+1}M^{2n+2}(t)=\int_{\R}h(t, x)g^2(t, x)dx\\
&\leq&\frac{2-c_1}{2}M^{2n}(t)\int_{\R}g^2(t, x)dx=\frac{2-c_1}{2}M^{2n}(t)\big(E(u)-2M^2(t)\big),
\end{eqnarray*}
which implies (\ref{3.5}). This completes the proof of Lemma 3.2.
\end{proof}

\begin{lemma3}
For $u\in H^s(\R),s>\frac{5}{2}$, if $\big\|u-\varphi_c\big\|_{H^{1}(\R)}<\varepsilon, \ \mbox{with} \ 0<\varepsilon<\big(3-2\sqrt{2}\big)a$, then
\begin{eqnarray*}
|E(u)-E(\varphi_c)|\leq 3a \varepsilon \quad \mbox{and} \quad |F(u)-F(\varphi_c)|
\lesssim  A(n,c,\|u\|_{H^s})\cdot \varepsilon,
 \end{eqnarray*}
where the constant $A(n,c,\|u\|_{H^s})>0$ depends only on $c>0$, integer $n>0$ and $\|u\|_{H^s}$.
\end{lemma3}
\begin{proof}
Under the assumption $0<\varepsilon<\big(3-2\sqrt{2}\big)a$, a direct computation yields
\begin{eqnarray}\label{3.28}
\big|E(u)-E(\varphi_c)\big|&=&\big|\big( \|u\|_{H^1}-\|\varphi_c\|_{H^1}\big)\big(
\|u\|_{H^1}+\|\varphi_c\|_{H^1}\big)\big|\nonumber\\
&\leq& \|u-\varphi_c\|_{H^1}\big(\|u-\varphi_c\|_{H^1}+2\|\varphi_c\|_{H^1}\big)\leq \varepsilon \big(\varepsilon+ 2\sqrt{2}a\big)
\leq 3a \varepsilon.
\end{eqnarray}
On the other hand, we have
\begin{eqnarray}\label{3.29}
&&|F(u)-F(\varphi_c)|\nonumber\\
&=&\Big|\int_{\R}\big(u^{2n+2}+\sum_{k=1}^{n}\frac{(-1)^{k+1}}{2k-1}C^k_{n+1}u^{2n-2k+2}u_x^{2k}+(-1)^n\frac{u^{2n+2}_x}{2n+1}\big)dx\nonumber\\
&&-\int_{\R}\big(\varphi_c^{2n+2}+\sum_{k=1}^{n}\frac{(-1)^{k+1}}{2k-1}C^k_{n+1}\varphi_c^{2n-2k+2}
(\partial_x\varphi_c)^{2k}+(-1)^n\frac{(\partial_x\varphi_c)^{2n+2}}{2n+1}\big)dx\Big|\nonumber\\
&\leq& \Big|\int_{\R}\big(u^{2n+2}+(n+1)u^{2n}u^2_x-\varphi_c^{2n+2}-(n+1)\varphi_c^{2n}(\partial_x\varphi_c)^2\big)dx\Big|+\frac{1}{2n+1}\nonumber\\
&&\times\Big|\int_{\R}\big(u^{2n+2}_x-(\partial_x\varphi_c)^{2n+2}\big)dx\Big|
+\sum_{k=2}^{n}\frac{C^k_{n+1}}{2k-1}\Big|\int_{\R}\big(u^{2n-2k+2}u_x^{2k}-\varphi_c^{2n-2k+2}
(\partial_x\varphi_c)^{2k}\big)dx\Big|\nonumber\\
&\leq& \Big|\int_{\R}\big(u^{2n}-\varphi_c^{2n}\big)\big(u^2+(n+1)u^2_x\big)dx\Big|
+\Big|\int_{\R}\varphi_c^{2n}\big((u^2-\varphi_c^2)+(n+1)(u^2_x-(\partial_x\varphi_c)^2)\big)dx\Big|\nonumber\\
&&+\frac{1}{2n+1}\Big|\int_{\R}\big(u_x^{n+1}+(\partial_x\varphi_c)^{n+1}\big)\big(u_x^{n+1}-(\partial_x\varphi_c)^{n+1}\big)dx\Big|
+\sum_{k=2}^{n}\frac{C^k_{n+1}}{2k-1}\Big|\int_{\R}u^{2n-2k+2}\nonumber\\
&&\times\big(u_x^{2k}-(\partial_x\varphi_c)^{2k}\big)dx\Big|+\sum_{k=2}^{n}\frac{C^k_{n+1}}{2k-1}
\Big|\int_{\R}(\partial_x\varphi_c)^{2k}\big(u^{2n-2k+2}-\varphi_c^{2n-2k+2}\big)dx\Big|\nonumber\\
&\triangleq& K_1+K_2+\cdot\cdot\cdot+K_5 \lesssim  A(n,c,\|u\|_{H^s})\cdot \varepsilon,
\end{eqnarray}
where the last inequality of (\ref{3.29}) will be proved in Appendix A.2. This completes the proof of Lemma 3.3.
\end{proof}
\begin{lemma3}
For $0<u\in H^s(\R),s>\frac{5}{2}$, let $M=\max_{x\in\R}\{u(x)\}.$ If
\begin{eqnarray*}
|E(u)-E(\varphi_c)|\leq 3a \varepsilon \quad \mbox{and} \quad |F(u)-F(\varphi_c)|
\lesssim  A(n,c,\|u\|_{H^s})\cdot \varepsilon,
 \end{eqnarray*}
with $0<\varepsilon<\big(3-2\sqrt{2}\big)a$ and the constant $A(n,c,\|u\|_{H^s})>0$ depends only on $c>0$, integer $n>0$
and $\|u\|_{H^s}$, then
\begin{eqnarray*}
\big|M-a\big|\lesssim  \sqrt{A(n,c,\|u\|_{H^s})\varepsilon}.
\end{eqnarray*}
\end{lemma3}
\begin{proof}
In view of (\ref{3.5}) in Lemma 3.2, we have
\begin{eqnarray}\label{3.31}
nM^{2n+2}-\frac{n+1}{2}M^{2n}E(u)+\frac{n+1}{2-c_1}F(u)\leq 0.
\end{eqnarray}
Consider the $(2n+2)$-th order polynomial $P(y)=ny^{2n+2}-\frac{n+1}{2}y^{2n}E(u)+\frac{n+1}{2-c_1}F(u).$  In the case
$E(u)=E(\varphi_c)=2a^2$ and $F(u)=F(\varphi_c)=\frac{1}{n+1}a^{2n+2} \big(1+\sum_{k=1}^{n+1}\frac{(-1)^{k+1}}{2k-1}C^k_{n+1}\big)$,
it takes the form
\begin{eqnarray}\label{3.32}
P_0(y)&=&ny^{2n+2} -\frac{n+1}{2}E(\varphi_c)y^{2n}+\frac{n+1}{2-c_1}F(\varphi_c)\nonumber\\
&=&ny^{2n+2}-(n+1)a^2y^{2n}+\frac{1}{2-c_1}\cdot a^{2n+2} \big(1+\sum\limits_{k=1}^{n+1}\frac{(-1)^{k+1}}{2k-1}C^k_{n+1}\big).
\end{eqnarray}
It follows from
\begin{eqnarray}\label{3.34}
\sum\limits_{k=1}^{n+1}(-1)^{k+1}C^k_{n+1}=1,
\end{eqnarray}
that the following identity holds,
\begin{eqnarray}\label{3.33}
2-c_1=1+\sum\limits_{k=1}^{n+1}\frac{(-1)^{k+1}}{2k-1}C^k_{n+1}.
\end{eqnarray}
Thus, the combination of (\ref{3.32}) and (\ref{3.33}) yields
\begin{eqnarray}\label{3.35}
P_0(y)&=&ny^{2n+2}-(n+1)a^2y^{2n}+a^{2n+2}\nonumber\\
&=&(y-a)^2\big(ny^{2n}+\sum_{k=1}^{2n-1}(2n+1-k)a^ky^{2n-k}+a^{2n}\big).
\end{eqnarray}
By (\ref{3.31}) and (\ref{3.32}), we obtain
\begin{eqnarray*}
P_0(M) \leq P_0(M)-P(M)=\frac{n+1}{2}M^{2n}\big(E(u)-E(\varphi_c)\big)-\frac{n+1}{2-c_1}\big(F(u)-F(\varphi_c)\big),
\end{eqnarray*}
which together with (\ref{3.35}) leads to
\begin{eqnarray}\label{3.36}
a^{2n}(M-a)^2 \leq \frac{n+1}{2}M^{2n}\big(E(u)-E(\varphi_c)\big)-\frac{n+1}{2-c_1}\big(F(u)-F(\varphi_c)\big).
\end{eqnarray}
By (\ref{3.6}) and the assumption of this lemma, we deduce that for $0<\varepsilon<\big(3-2\sqrt{2}\big)a$,
\begin{eqnarray}\label{3.37}
0<M^2 \leq \frac{E(u)}{2}\leq \frac{2a^2+3a\varepsilon}{2}<\frac{3a^2}{2}.
\end{eqnarray}
Hence, combining (\ref{3.36}) with (\ref{3.37}), we find that
\begin{eqnarray*}
a^n\big|M-a\big|\lesssim \sqrt{\frac{(n+1)3^{n+1}}{2^{n+1}}a^{2n+1}\varepsilon+\frac{n+1}{2}A(n,c,\|u\|_{H^s})\varepsilon}.
\end{eqnarray*}
Thus, it follows from the above inequality and the relation $c=\big(1-\sum_{k=1}^n\frac{(-1)^{k+1}}{2k+1}C^k_n\big)a^{2n}$ that there is
a constant, still denoted by $A(n,c,\|u\|_{H^s})$ for simplicity, such that
\begin{eqnarray*}
\big|M-a\big|\lesssim  \sqrt{A(n,c,\|u\|_{H^s})\varepsilon}.
\end{eqnarray*}
This completes the proof of Lemma 3.4.
\end{proof}

We are now in a position to prove Theorem 3.1.\\

$\textbf{Proof \ of \ Theorem \ 3.1.}$  Let $u\in C([0,T);H^{s}(\R)), s>\frac{5}{2}$ be the solution of
equation (\ref{1.1}) with the initial data $u_0\in H^s(\R)$. Since $E(u)$ and $F(u)$ are both conserved densities of equation
(\ref{1.1}), we thus have
\begin{eqnarray}\label{3.38}
E(u(t,\cdot))=E(u_0) \quad \mbox{and} \quad  F(u(t,\cdot))=F(u_0),
\quad \forall \ t\in [0,T).
\end{eqnarray}
Since $\|u(0,\cdot)-\varphi_c\|_{H^{s}}<\varepsilon, \ \mbox{with} \ 0<\varepsilon<\big(3-2\sqrt{2}\big)a$, and $0\neq(1-\partial^2_x)u_0\geq 0$,
in view of (\ref{3.38}) and Lemma 3.3, the hypotheses of Lemma 3.4 are satisfied for $u(t,\cdot)$ with a positive constant $A(n,c,\|u_0\|_{H^s})$ to be chosen depending only on wave speed $c>0$, integer $n>0$ and $\|u_0\|_{H^s}$. Hence
\begin{eqnarray}\label{3.39}
\big|u(t,\xi(t))-a\big|\lesssim  \sqrt{A(n,c,\|u_0\|_{H^s})\varepsilon}, \quad \forall \ t\ \in[0,T),
\end{eqnarray}
where $\xi(t) \in {\R}$ is the maximum point of function $u(t,\cdot)$. Combining (\ref{3.4}) with (\ref{3.38}), we find
\begin{eqnarray*}
\big\|u-\varphi_c(\cdot-\xi(t))\big\|^2_{H^1}=E(u_0)-E(\varphi_c)
-4 a\big(u(t,\xi(t))-a\big).
\end{eqnarray*}
Therefore, it follows (\ref{3.39}) and Lemma 3.3 that for $t\in [0,T)$,
\begin{eqnarray*}
\big\|u-\varphi_c(\cdot-\xi(t))\big\|_{H^1}\leq
\sqrt{\big|E(u_0)-E(\varphi_c)\big|+4 a\big|u(t,\xi(t))-a\big|}\lesssim  \sqrt{3a\varepsilon+4a\sqrt{A(c,\|u_0\|_{H^s})\varepsilon}}.
\end{eqnarray*}
This completes the proof of Theorem 3.1.

\section{Appendix}

In this section, we provide the details of the proofs to the last third equality (\ref{3.8}), last inequality of (\ref{3.29}) and the conservation law $F(u)$.
\\
$\textbf{A.1.\ Proof\ of\ the\ last\ third\ equality\ of\ (\ref{3.8}).}$ In Lemma 3.2, we define the function $h(x)$ as
\begin{eqnarray*}
h(x)=\left\{\begin{array}{ll}
\big(u^{2n}+\sum\limits_{k=1}^{2n-1}c_ku^{2n-k}u_x^k+(-1)^n\frac{1}{2n+1}u^{2n}_x\big)(x), & x<\xi,\\
\big(u^{2n}+\sum\limits_{k=1}^{2n-1}d_ku^{2n-k}u_x^k+(-1)^n\frac{1}{2n+1}u^{2n}_x\big)(x),& x>\xi.
\end{array}\right.
\end{eqnarray*}
To determine the coefficients $c_k,d_k,k=1,2,\cdot\cdot\cdot,2n-1$, as in (\ref{3.7}), a useful observation is that $h(x)$ satisfies the last third equality of (\ref{3.8}) .
Hence, we do calculation as follows
\begin{eqnarray}\label{4.1}
&&(u^2-2uu_x+u^2_x)\big(u^{2n}+\sum\limits_{k=1}^{2n-1}c_ku^{2n-k}u_x^k+(-1)^n\frac{1}{2n+1}u^{2n}_x\big)\nonumber\\
&=&u^{2n+2}+(c_1-2)u^{2n+1}u_x+(c_2-2c_1+1)u^{2n}u^2_x+\sum_{k=3}^{2n-1}(c_k-2c_{k-1}+c_{k-2})u^{2n-k+2}u^k_x\nonumber\\
&&+\big(\frac{(-1)^n}{2n+1}-2c_{2n-1}+c_{2n-2}\big)u^2u^{2n}_x+\big(-\frac{2\cdot(-1)^n}{2n+1}+c_{2n-1}\big)uu^{2n+1}_x
+\frac{(-1)^n}{2n+1}u^{2n+2}_x,\nonumber\\
\end{eqnarray}
and
\begin{eqnarray}\label{4.2}
&&(u^2+2uu_x+u^2_x)\big(u^{2n}+\sum\limits_{k=1}^{2n-1}d_ku^{2n-k}u_x^k+(-1)^n\frac{1}{2n+1}u^{2n}_x\big)\nonumber\\
&=&u^{2n+2}+(d_1+2)u^{2n+1}u_x+(d_2+2d_1+1)u^{2n}u^2_x+\sum_{k=3}^{2n-1}(d_k+2d_{k-1}+d_{k-2})u^{2n-k+2}u^k_x\nonumber\\
&&+\big(\frac{(-1)^n}{2n+1}+2d_{2n-1}+d_{2n-2}\big)u^2u^{2n}_x+\big(\frac{2\cdot(-1)^n}{2n+1}+d_{2n-1}\big)uu^{2n+1}_x
+\frac{(-1)^n}{2n+1}u^{2n+2}_x.\nonumber\\
\end{eqnarray}
For our purpose, combining the conservation law $F(u)$ with (\ref{4.1})-(\ref{4.2}), we thus obtain
\begin{equation}\label{4.3}
\left\{\begin{array}{llll}-\frac{2\cdot(-1)^n}{2n+1}+c_{2n-1}=0,\\
\frac{(-1)^n}{2n+1}-2c_{2n-1}+c_{2n-2}=\frac{(-1)^{n+1}}{2n-1}C_{n+1}^n,\\
c_k-2c_{k-1}+c_{k-2}=0,& k=2j+1,\ (j=1,2,\cdot\cdot\cdot,n-1),\\
c_k-2c_{k-1}+c_{k-2}=\frac{(-1)^{j+1}}{2j-1}C_{n+1}^j,& k=2j,\ (j=2,3,\cdot\cdot\cdot,n-1),\\
c_{2}-2c_{1}+1=C_n^1+C_n^0,
\end{array}\right.
\end{equation}
and
\begin{equation}\label{4.4}
\left\{\begin{array}{llll}\frac{2\cdot(-1)^n}{2n+1}+d_{2n-1}=0,\\
\frac{(-1)^n}{2n+1}+2d_{2n-1}+d_{2n-2}=\frac{(-1)^{n+1}}{2n-1}C_{n+1}^n,\\
d_k+2d_{k-1}+d_{k-2}=0,& k=2j+1,\ (j=1,2,\cdot\cdot\cdot,n-1),\\
d_k+2d_{k-1}+d_{k-2}=\frac{(-1)^{j+1}}{2j-1}C_{n+1}^j,& k=2j,\ (j=2,3,\cdot\cdot\cdot,n-1),\\
d_{2}+2d_{1}+1=C_n^1+C_n^0.
\end{array}\right.
\end{equation}
By induction, (\ref{3.7}) follows from the relations (\ref{4.3})-(\ref{4.4}). Hence, the last third equality in (\ref{3.8}) also holds.\\
$\textbf{A.2.\ Proof\ of\ the\ last\ inequality\ of\ (\ref{3.29}).}$ To complete the proof of Lemma 3.3, it remains to prove the following inequality:
\begin{eqnarray}\label{4.5}
K_1+K_2+\cdot\cdot\cdot+K_5 \lesssim  A(n,c,\|u\|_{H^s})\cdot \varepsilon.
\end{eqnarray}
Using the relation (\ref{3.6}), we have for all $v\in H^1(\R)$,
\begin{eqnarray}\label{4.6}
\sup_{x\in \R}|v(x)|\leq \frac{\sqrt{E(v)}}{\sqrt{2}}
\leq\frac{\|v\|_{H^1}}{\sqrt{2}}.
\end{eqnarray}
For the term $K_1$, by (\ref{3.28}) and (\ref{4.6}), we have
\begin{eqnarray*}
K_1&\leq&(n+1)\int_{\R}\big|u-\varphi_c\big|\cdot
\big|u^{2n-1}+u^{2n-2}\varphi_c+\cdot\cdot\cdot
+u\varphi_c^{2n-2}+\varphi_c^{2n-1}\big|\cdot\big(u^2+u^2_x\big)dx\\
&\leq&(n+1) E(u) \|u-\varphi_c\|_{L^\infty}\big(\|u\|_{L^\infty}^{2n-1}+\|u\|_{L^\infty}^{2n-2}\|\varphi_c\|_{L^\infty}
+\cdot\cdot\cdot+\|u\|_{L^\infty}\|\varphi_c\|^{2n-2}_{L^\infty}\\
&&+\|\varphi_c\|^{2n-1}_{L^\infty}\big)\\
&\leq&\frac{n+1}{2^n} \big(E(\varphi_c)+3a\varepsilon\big) \|u-\varphi_c\|_{H^1}\big(\|u\|_{H^1}^{2n-1}+\|u\|_{H^1}^{2n-2}\|\varphi_c\|_{H^1}
+\cdot\cdot\cdot+\|u\|_{H^1}\|\varphi_c\|^{2n-2}_{H^1}\\
&&+\|\varphi_c\|^{2n-1}_{H^1}\big)\\
&\leq& \frac{n+1}{2^n}\big(E(\varphi_c)+3a \varepsilon\big) \big(\|u-\varphi_c\|_{H^1}+2\|\varphi_c\|_{H^1}\big)^{2n-1}\cdot\|u-\varphi_c\|_{H^1} \\
&\leq &  \frac{n+1}{2^n} (2a^2+3a \varepsilon) (2\sqrt{2}a+\varepsilon)^{2n-1}\cdot\varepsilon.
\end{eqnarray*}
Similarly, for $K_2$, we have
\begin{eqnarray*}
K_2&\leq&\|\varphi_c\|_{L^\infty}^{2n}\Big|\int_{\R}\big((u-\varphi_c)^2+(n+1)(u_x-\partial_x\varphi_c)^2
+2\varphi_c(u-\varphi_c)\\
&&+2(n+1)\partial_x\varphi_c(u_x-\partial_x\varphi_c)\big)dx\Big|\\
&\leq&\big(\frac{\|\varphi_c\|_{H^1}}{\sqrt{2}}\big)^{2n}
\big((n+1)\|u-\varphi_c\|^2_{H^1}+2(n+1)\|\varphi_c\|_{H^1}\|u-\varphi_c\|_{H^1}\big)\\
&\leq &  (n+1)a^{2n} (2\sqrt{2}a+ \varepsilon) \cdot \varepsilon.
\end{eqnarray*}
For the term $K_3$, by the H\"{o}lder inequality, we obtain
\begin{eqnarray*}
K_3&=&\frac{1}{2n+1}\Big|\int_{\R}\big(u_x^{n+1}+(\partial_x\varphi_c)^{n+1}\big)\big(u_x^n+
u_x^{n-1}(\partial_x\varphi_c)+\cdot\cdot\cdot+u_x(\partial_x\varphi_c)^{n-1}\\
&&+(\partial_x\varphi_c)^{n}
\big)\cdot\big(u_x-\partial_x\varphi_c\big)dx\Big|\\
&\leq&\frac{1}{2n+1} \Big(\int_{\R}\big(u_x^{n+1}+(\partial_x\varphi_c)^{n+1}\big)^2\big(u_x^n+
u_x^{n-1}(\partial_x\varphi_c)+\cdot\cdot\cdot+u_x(\partial_x\varphi_c)^{n-1}\\
&&+(\partial_x\varphi_c)^{n}
\big)^2dx\Big)^{\frac{1}{2}}
 \Big(\int_{\R}\big(u_x-\partial_x\varphi_c\big)^2dx\Big)^{\frac{1}{2}}\\
&\triangleq& \frac{1}{2n+1}\sqrt{K'_{3}} \|u-\varphi_c\|_{H^1}.
\end{eqnarray*}
For $K'_{3}$, a direct use of Young's inequality yields,
\begin{eqnarray*}
K'_{3}&=&\int_{\R}\big(u_x^{4n+2}+2u_x^{4n+1}\partial_x\varphi_c+3u_x^{4n}(\partial_x\varphi_c)^2
+\cdot\cdot\cdot+(2n+1)u_x^{2n+2}(\partial_x\varphi_c)^{2n}\\
&&+(2n+2)u_x^{2n+1}(\partial_x\varphi_c)^{2n+1}
+(2n+1)u_x^{2n}(\partial_x\varphi_c)^{2n+2}\\
&&+\cdot\cdot\cdot+3u_x^2(\partial_x\varphi_c)^{4n}+2u_x(\partial_x\varphi_c)^{4n+1}+(\partial_x\varphi_c)^{4n+2}\big)dx\\
&\leq& 2(n+1)^2\big(\int_{\R} u_x^{4n+2} dx+\int_{\R} (\partial_x\varphi_c)^{4n+2} dx\big).
\end{eqnarray*}
Since $u\in H^s(\R)\subset H^2(\R),s>\frac{5}{2}$, by the Gagliardo-Nirenberg inequality, we have
\begin{eqnarray*}
\|u_x\|^{4n+2}_{L^{4n+2}}\leq C \|u\|^{n+1}_{L^{2}}\|u_{xx}\|_{L^2}^{3n+1},
\end{eqnarray*}
with the constant $C>0$ independent of $u.$ Hence it follows from $\|\partial_x\varphi_c\|^{4n+2}_{L^{4n+2}}=\frac{a^{4n+2}}{2n+1}$ that there holds
\begin{eqnarray*}
K'_{3}\lesssim A^2(n,a,\|u\|_{H^s}),
\end{eqnarray*}
where the constant $A(n,a,\|u\|_{H^s})>0$ depends only on $a>0$, integer $n>0$, and the norm $\|u\|_{H^s}$.
Since $ c=\big(1-\sum_{k=1}^n\frac{(-1)^{k+1}}{2k+1}C^k_n\big)a^{2n}$, then we have
\begin{eqnarray*}
K_3\lesssim  A(n,c,\|u\|_{H^s})\cdot \varepsilon.
\end{eqnarray*}
Applying the similar method to treat $K_4$, by the fact that
$\|\partial_x\varphi_c\|^{4k-2}_{L^{4k-2}}=\frac{a^{4k-2}}{2k-1}$, we obtain
\begin{eqnarray*}
K_4&\leq& \sum_{k=2}^{n}\frac{C^k_{n+1}}{2k-1}\|u\|_{L^\infty}^{2n-2k+2}
\Big(\int_{\R}\big(u_x^k+(\partial_x\varphi_c)^k\big)^2
\big(u^{k-1}_x+u^{k-2}_x\partial_x\varphi_c+\cdot\cdot\cdot+u_x(\partial_x\varphi_c)^{k-2}\\
&&+(\partial_x\varphi_c)^{k-1}
\big)^2dx\Big)^{\frac{1}{2}}
 \Big(\int_{\R}\big(u_x-\partial_x\varphi_c\big)^2dx\Big)^{\frac{1}{2}}\\
&=& \sum_{k=2}^{n}\frac{C^k_{n+1}}{2k-1}\|u\|_{L^\infty}^{2n-2k+2}
\Big(\int_{\R}\big(u_x^{4k-2}+2u^{4k-3}_x\partial_x\varphi_c+\cdot\cdot\cdot+
(2k-1)u^{2k}_x(\partial_x\varphi_c)^{2k-2}\\
&&+2ku^{2k-1}_x(\partial_x\varphi_c)^{2k-1}
+(2k-1)u^{2k-2}_x(\partial_x\varphi_c)^{2k}+\cdot\cdot\cdot+2 u_x(\partial_x\varphi_c)^{4k-3}\\
&&+(\partial_x\varphi_c)^{4k-2}\big)dx\Big)^{\frac{1}{2}}\Big(\int_{\R}\big(u_x-\partial_x\varphi_c\big)^2dx\Big)^{\frac{1}{2}}\\
&\leq& \sum_{k=2}^{n}\frac{C^k_{n+1}}{2k-1}\|u\|_{L^\infty}^{2n-2k+2}\cdot \sqrt{2}k\big(\int_{\R} u_x^{4k-2} dx+\int_{\R} (\partial_x\varphi_c)^{4k-2} dx\big)^{\frac{1}{2}}\cdot \|u-\varphi_c\|_{H^1}\\
&\lesssim & \|u\|_{H^1}^{2n-2k+2} \big(\|u\|^{k}_{L^2}\|u_{xx}\|_{L^2}^{3k-2}+\|\partial_x\varphi_c\|^{4k-2}_{L^{4k-2}}\big)^{\frac{1}{2}}\cdot \|u-\varphi_c\|_{H^1}\lesssim A(n,c,\|u\|_{H^s})\cdot \varepsilon.
\end{eqnarray*}
For $K_5$, a direct use of the H\"{o}lder inequality gives rise to,
\begin{eqnarray*}
K_5&\leq& \sum_{k=2}^{n}\frac{C^k_{n+1}}{2k-1} \big(\|u\|_{L^\infty}^{n-k+1} +\|\varphi_c\|_{L^\infty}^{n-k+1}\big)
\big(\|u\|^{n-k}_{L^\infty} +\|u\|^{n-k-1}_{L^\infty}\|\varphi_c\|_{L^\infty}+\cdot\cdot\cdot\\
&&+\|u\|_{L^\infty}\|\varphi_c\|^{n-k-1}_{L^\infty}+\|\varphi_c\|^{n-k}_{L^\infty}\big)
\Big(\int_{\R}\big(\partial_x\varphi_c\big)^{4k}dx\Big)^{\frac{1}{2}}
 \Big(\int_{\R}\big(u-\varphi_c\big)^2dx\Big)^{\frac{1}{2}}\\
 &\leq &\sum_{k=2}^{n}\frac{C^k_{n+1}}{2k-1}\big(\|u-\varphi_c\|_{H^1}+2\|\varphi_c\|_{H^1}\big)^{2n-2k+1}\|\partial_x\varphi_c\|^{2k}_{L^{4k}}\cdot \|u-\varphi_c\|_{H^1}\\
 &\leq &\sum_{k=2}^{n}\frac{C^k_{n+1}}{2k-1}\cdot\frac{a^{4k}}{2k}(2\sqrt{2}a+\varepsilon)^{2n-2k+1}\cdot \varepsilon,
\end{eqnarray*}
where we used the fact that $\|\partial_x\varphi_c\|^{4k}_{L^{4k}}=\frac{a^{4k}}{2k}$. Therefore,
gathering the estimations $K_1$-$K_5$, we arrive at the desired result (\ref{4.5}).\\
$\textbf{A.3.\ Proof\ of\ the\ conservation\ law\ $F(u)$.}$ 
Setting $v(x,t)=\int^x_{-\infty}u_t(z,t)dz$, it follows form integration by parts that 
\begin{eqnarray}\label{2.24}
\frac{d}{dt}\int_{\R}u^{2n+2}dx=(2n+2)\int_{\R}u^{2n+1}v_xdx
=-(2n+2)\int_{\R}(2n+1)u^{2n}u_xvdx,
\end{eqnarray}
\begin{eqnarray}\label{2.25}
\frac{d}{dt}\int_{\R}\frac{(-1)^n}{2n+1}u^{2n+2}_xdx
=(2n+2)\int_{\R}(-1)^n\big(2nu_x^{2n-1}u_{xx}^2+u_x^{2n}u_{xxx}\big)vdx,
\end{eqnarray}
and
\begin{eqnarray}\label{2.26}
&&\frac{d}{dt}\int_{\R}(n+1)u^{2n}u^2_xdx\nonumber\\
&=&(n+1)\int_{\R}\big(2nu^{2n-1}u_x^2 u_t+2u^{2n}u_xu_{tx}\big)dx
=(2n+2)\int_{\R}(nu^{2n-1}u_x^2v_x+u^{2n}u_xv_{xx})dx\nonumber\\
&=& (2n+2)\int_{\R}\big(n(2n-1)u^{2n-2}u^3_x+4nu^{2n-1}u_xu_{xx}+u^{2n}u_{xxx}\big)vdx.
\end{eqnarray}
In a similar manner,
\begin{eqnarray}\label{2.27}
&&\frac{d}{dt}\int_{\R}\sum_{k=2}^{n}\frac{(-1)^{k+1}}{2k-1}C^k_{n+1}u^{2n-2k+2}u_x^{2k}dx\nonumber\\
&=&\sum_{k=2}^{n}\frac{(-1)^{k+1}}{2k-1}C^k_{n+1} \int_{\R}\big(2(n-k+1)u^{2n-2k+1}u_x^{2k}u_t+2ku^{2(n-k+1)}u^{2k-1}_xu_{tx}\big)dx\nonumber\\
&=&\sum_{k=2}^{n}\frac{(-1)^{k+1}}{2k-1}C^k_{n+1} \int_{\R}\big(2(n-k+1)u^{2n-2k+1}u_x^{2k}v_x+2ku^{2(n-k+1)}u^{2k-1}_xv_{xx}\big)dx\nonumber\\
&=&\sum_{k=2}^{n}\frac{(-1)^{k+1}}{2k-1}C^k_{n+1}\int_{\R}\Big[-2(n-k+1)\big((2n-2k+1)u^{2n-2k}u^{2k+1}_x\nonumber\\
&&+2ku^{2n-2k+1}u^{2k-1}_xu_{xx}\big)v dx+\big(4k(n-k+1)(2n-2k+1)u^{2n-2k}u^{2k+1}_x \nonumber\\
&&+2(n-k+1)(8k^2-2k)u^{2n-2k+1}u^{2k-1}_xu_{xx}+2k(2k-1)(2k-2)u^{2n-2k+2}u^{2k-3}_xu^2_{xx}\nonumber\\
&&+2k(2k-1)u^{2n-2k+2}u^{2k-2}_xu_{xxx}\big)vdx\Big]\nonumber\\
&=&\sum_{k=2}^{n}(-1)^{k+1}C^k_{n+1}\int_{\R}\big(2(n-k+1)(2n-2k+1)u^{2n-2k}u^{2k+1}_x\nonumber\\
&&+2(n-k+1)4ku^{2n-2k+1}u^{2k-1}_xu_{xx}+2k(2k-2)u^{2n-2k+2}u^{2k-3}_xu^2_{xx}\nonumber\\
&&+2ku^{2n-2k+2}u^{2k-2}_xu_{xxx})vdx\big)\triangleq J_1+J_2+J_3+J_4.
\end{eqnarray}
For $J_1$-$J_4$, a direct computation yields
\begin{eqnarray*}
J_1&=&\sum_{k=2}^{n}(-1)^{k+1}C^k_{n+1}\int_{\R}2(n-k+1)(2n-2k+1)u^{2n-2k}u^{2k+1}_xdx\\
&=&\sum_{k=2}^{n}(-1)^{k+1} \frac{(n+1)!}{k!(n-k+1)!}2(n-k+1)(2n-2k+1) \int_{\R} u^{2n-2k}u^{2k+1}_xdx\\
&=&(2n+2)\sum_{k=2}^{n}(-1)^{k+1} \frac{n!}{k!(n-k)!}(2n-2k+1) \int_{\R} u^{2n-2k}u^{2k+1}_xdx\\
&=&(-1)^{n+1}(2n+2)\int_{\R}u^{2n+1}_xdx+(2n+2)\sum_{k=2}^{n-1}(-1)^{k+1} C_n^k(2n-2k+1) \int_{\R} u^{2n-2k}u^{2k+1}_xdx\\
&=&(-1)^{n+1}(2n+2)\int_{\R}u^{2n+1}_xdx+(2n+2)\sum_{k=2}^{n-1}(-1)^{k+1} \big(C_n^k+2nC_{n-1}^k\big)\int_{\R} u^{2n-2k}u^{2k+1}_xdx,
\end{eqnarray*}
\begin{eqnarray*}
J_2&=&\sum_{k=1}^{n-1}(-1)^{k}C^{k+1}_{n+1}\int_{\R}2(n-k)4(k+1)u^{2n-2k-1}u^{2k+1}_xu_{xx}dx\\
&=&\sum_{k=1}^{n-1}(-1)^{k} \frac{(n+1)!}{(k+1)!(n-k)!}2(n-k)4(k+1)\int_{\R}u^{2n-2k-1}u^{2k+1}_xu_{xx}dx\\
&=&(2n+2)\sum_{k=1}^{n-1}(-1)^{k} \frac{(n-1)!}{k!(n-k-1)!}4n \int_{\R}u^{2n-2k-1}u^{2k+1}_xu_{xx}dx\\
&=&(2n+2)\sum_{k=1}^{n-1}(-1)^{k} 4n C_{n-1}^k \int_{\R}u^{2n-2k-1}u^{2k+1}_xu_{xx}dx,
\end{eqnarray*}
\begin{eqnarray*}
J_3&=&\sum_{k=0}^{n-2}(-1)^{k+1}C^{k+2}_{n+1}\int_{\R}4(k+2)(k+1)u^{2n-2k-2}u^{2k+1}_xu^2_{xx}dx\\
&=&\sum_{k=0}^{n-2}(-1)^{k+1} \frac{(n+1)!}{(k+2)!(n-k-1)!}4(k+2)(k+1)\int_{\R}u^{2n-2k-2}u^{2k+1}_xu^2_{xx}dx\\
&=&(2n+2)\sum_{k=0}^{n-2}(-1)^{k+1} \frac{(n-1)!}{k!(n-k-1)!}2n \int_{\R}u^{2n-2k-2}u^{2k+1}_xu^2_{xx}dx\\
&=&(2n+2)\sum_{k=0}^{n-2}(-1)^{k+1} 2nC_{n-1}^k \int_{\R}u^{2n-2k-2}u^{2k+1}_xu^2_{xx}dx,
\end{eqnarray*}
and
\begin{eqnarray*}
J_4&=&\sum_{k=1}^{n-1}(-1)^{k}C^{k+1}_{n+1}\int_{\R}2(k+1)u^{2n-2k}u^{2k}_xu_{xxx}dx\\
&=&\sum_{k=1}^{n-1}(-1)^{k} \frac{(n+1)!}{(k+1)!(n-k)!}2(k+1)\int_{\R}u^{2n-2k}u^{2k}_xu_{xxx}dx\\
&=&(2n+2)\sum_{k=1}^{n-1}(-1)^{k} \frac{n!}{k!(n-k)!} \int_{\R}u^{2n-2k}u^{2k}_xu_{xxx}dx\\
&=&(2n+2)\sum_{k=1}^{n-1}(-1)^{k} C_n^k \int_{\R}u^{2n-2k}u^{2k}_xu_{xxx}dx.
\end{eqnarray*}
Plugging $J_1$-$J_4$ into (\ref{2.27}), we deduce that
\begin{eqnarray}\label{2.28}
&&\frac{d}{dt}\int_{\R}\sum_{k=2}^{n}\frac{(-1)^{k+1}}{2k-1}C^k_{n+1}u^{2n-2k+2}u_x^{2k}dx\nonumber\\
&=&(2n+2)\int_{\R}\Big[(-1)^{n+1}u^{2n+1}_x+\sum_{k=2}^{n-1}(-1)^{k+1} \big(C_n^k+2nC_{n-1}^k\big) u^{2n-2k}u^{2k+1}_x\nonumber\\
&&+\sum_{k=1}^{n-1}(-1)^{k} 4n C_{n-1}^k u^{2n-2k-1}u^{2k+1}_xu_{xx}
+\sum_{k=0}^{n-2}(-1)^{k+1} 2nC_{n-1}^k u^{2n-2k-2}u^{2k+1}_xu^2_{xx}\nonumber\\
&&+\sum_{k=1}^{n-1}(-1)^{k} C_n^k u^{2n-2k}u^{2k}_xu_{xxx}\Big]dx.
\end{eqnarray}
Combining (\ref{2.24})-(\ref{2.26}) with (\ref{2.28}), we have
\begin{eqnarray*}\frac{dF(u)}{dt}&=&\frac{d}{dt}\int_{\R}\Big(u^{2n+2}+\sum_{k=1}^{n}\frac{(-1)^{k+1}}{2k-1}C^k_{n+1}u^{2n-2k+2}u_x^{2k}+(-1)^n\frac{u^{2n+2}_x}{2n+1}\Big)dx\\
&=&-(2n+2)\int_{\R}\Big[
\sum_{k=0}^{n-1}(-1)^k\big(C^k_n+2nC_{n-1}^k\big)u^{2n-2k}u_x^{2k+1}+(-1)^nu_x^{2n+1}\\
&&-\sum_{k=0}^{n-1}(-1)^k4nC_{n-1}^ku^{2n-2k-1}u_x^{2k+1}u_{xx}+
\sum_{k=0}^{n-1}(-1)^k2nC_{n-1}^ku^{2n-2k-2}u_x^{2k+1}u^2_{xx}\\
&&-\sum_{k=0}^{n}(-1)^kC_{n}^ku^{2n-2k}u_x^{2k}u_{xxx} \big)v dx,
\end{eqnarray*}
which along with (\ref{2.0}) yields,
\begin{eqnarray*}\frac{dF(u)}{dt}=(2n+2)\int_{\R}(u_t-u_{txx})vdx=\int_{\R}(2n+2)(vv_x-vv_{xxx})dx=0.\end{eqnarray*}

\bigskip
\noindent\textbf{Acknowledgments} \
Z. Guo was  partially supported by ARC DP170101060. X.C. Liu was partially supported by NSFC under Grant 11722111 and Grant 11631007. X.X. Liu was partially supported by the Fundamental Research Funds for the Central Universities (No.2018QNA34 and No.2017XKZD11).
C. Qu was partially supported by NSFC under Grant 11471174 and Grant 11631007. The authors thank the anonymous
referee for helpful suggestions and comments.

\end{document}